\documentclass[11pt,a4paper]{article}
\usepackage[margin=29mm]{geometry}
\usepackage{amsmath,amssymb,amsthm,mathtools}
\usepackage{microtype}
\usepackage[hidelinks]{hyperref}
\allowdisplaybreaks
\numberwithin{equation}{section}
\newtheorem{theorem}{Theorem}[section]
\newtheorem{lemma}{Lemma}[section]
\newtheorem{proposition}[lemma]{Proposition}

\theoremstyle{definition}
\newtheorem{assumption}{Structural condition}[section]
\theoremstyle{remark}
\newtheorem{remark}[lemma]{Remark}
\newcommand{\R}{\mathbb R}
\newcommand{\curl}{\operatorname{curl}}
\newcommand{\TT}{\mathbb T^2}
\newcommand{\dd}{\,\mathrm dx}
\newcommand{\diver}{\operatorname{div}}
\newcommand{\Dt}{D_t}

\title{Global Strong Solutions and Large-Time Behavior
  for the Two-Dimensional Compressible Flow  with General Bulk Viscosity\thanks{This work is partially supported by the National Natural Science Foundation of
China (12471227, 12494541).}}

\author{Xinyu  F{\small AN}$^{a}$,  Song J{\small IANG}$^{b}$, Jing L{\small I}$^{c,d}$, Boqiang L{\small \"U}$^{d} $  \thanks{email:  fanxinyu17@mails.ucas.ac.cn(X. Fan),  jiang@iapcm.ac.cn(S. Jiang), ajingli@gmail.com (J. Li), lvbq86@ncu.edu.cn (B.  L\"u) }  \\
{\normalsize a. School of Mathematics  Sciences, } \\
{\normalsize  Nankai University, Tianjin 300110, P. R. China;}
 \\ {\normalsize b. Institute of Applied Physics and Computational Mathematics, }\\\ {\normalsize Beijing 100088, P. R. China;}\\
{\normalsize c. Institute of Applied Mathematics \& State Key Laboratory of Mathematical Sciences,  AMSS,} \\ {\normalsize \&   Hua Loo-Keng Key Laboratory of Mathematics,}\\
{\normalsize  Chinese Academy of Sciences,    Beijing 100190, P. R. China};\\
{\normalsize d. School of Mathematics and Computer Sciences,}\\  {\normalsize  Nanchang University, Nanchang 330031, P. R. China.}
}
\date{}

\begin{document}
\maketitle

\begin{abstract}
The global existence of a unique strong solution to the two-dimensional
barotropic compressible Navier--Stokes equations with general
density-dependent bulk viscosity is established for large initial data
allowing vacuum. Both the periodic problem and the Cauchy problem with
far-field vacuum are considered under a common logarithmic derivative
bound and growth condition on the viscosity primitive. A time-independent
density upper bound and uniform low-order estimates are obtained in both
domains. The admissible coefficients include every superlinear power law
and near-linear laws with logarithmic gains. Exponential relaxation to the
positive equilibrium is proved in the periodic setting. In the whole
space, uniform relaxation toward vacuum and algebraic decay of the
pressure, velocity gradient, material acceleration, and relative
kinetic energy are established.
\end{abstract}


\section{Introduction and main results}\label{sec:introduction}

We study global strong solutions and large-time behavior for
two-dimensional compressible flows with general density-dependent bulk
viscosity and large initial data allowing vacuum. The periodic and
whole-space problems use the same equations, viscosity primitive, and
growth condition; their initial data are specified separately below.

Let $\Omega=\TT:=\mathbb R^2/\mathbb Z^2$ or $\Omega=\R^2$.
We consider
\begin{equation}\label{eq:system}
\begin{cases}
 \rho_t+\diver(\rho u)=0,\\
 (\rho u)_t+\diver(\rho u\otimes u)+\nabla P
 =\mu\Delta u+\nabla((\mu+\lambda(\rho))\diver u),
\end{cases}
\end{equation}
where $\rho$ and $u$ denote the density and velocity, and $\mu>0$ and
$\gamma>1$ are fixed. We use
\[
 P=\rho^\gamma,\qquad \Dt=\partial_t+u\cdot\nabla,\qquad \dot u=\Dt u.
\]
Our curl convention is
\[
 \omega=\curl u=\partial_2u_1-\partial_1u_2,\qquad \nabla^\perp=(\partial_2,-\partial_1).
\]
The initial density and momentum are prescribed by
\begin{equation}\label{eq:initial-traces}
 \rho|_{t=0}=\rho_0,\qquad (\rho u)|_{t=0}=\rho_0u_0.
\end{equation}
Spatial integrals and norms are taken over the domain under discussion.
On the unit torus, $\bar f=\int_{\TT}f\dd$.

The following structural condition is common to both domains.
\begin{assumption}[Bulk viscosity]\label{per:ass:lambda}
The coefficient satisfies
\begin{equation}\label{eq:lambda-regularity}
 \lambda\in C^1([0,\infty))\cap C^2((0,\infty)),\qquad
 \lambda\ge0,\qquad
 r|\lambda'(r)|\le C_\lambda(2\mu+\lambda(r)),\quad r>0.
\end{equation}
Define the viscosity primitive by
\begin{equation}\label{eq:theta}
 \theta(r)=\int_1^r\frac{2\mu+\lambda(z)}z\,\mathrm dz
 =2\mu\log r+\int_1^r\frac{\lambda(z)}z\,\mathrm dz,
 \qquad r>0.
\end{equation}
For some $\kappa\in(1/2,1]$, $c_\lambda>0$, and $r_\lambda\ge2$,
\begin{equation}\label{per:eq:primitive-growth}
 \theta(r)\ge c_\lambda r[\log(e+r)]^\kappa,
 \qquad r\ge r_\lambda.
\end{equation}
\end{assumption}
We enlarge $C_\lambda$ so that the derivative bound also holds with
$1+\lambda$ on the right. The coefficient need not be monotone, and
$\lambda(0)$ may be positive. Then $\theta'(r)=(2\mu+\lambda(r))/r>0$.

For the power law $\lambda(\rho)=\rho^\beta$, Vaigant and
Kazhikhov~\cite{VaigantKazhikhov1995} established global strong
solvability for large periodic data bounded away from vacuum when
$\beta>3$. Jiu, Wang, and Xin~\cite{JiuWangXin2014} allowed initial
vacuum in the classical periodic problem under the same restriction on
$\beta$. Huang and Li~\cite[Theorems~1.1--1.2]{HuangLi2016} proved
global existence of strong and weak solutions for $\beta>4/3$ and
$\gamma>1$, including vacuum and without a smallness assumption on the
initial data. Fan, Li, and Wang~\cite[Theorem~1.2]{FanLiWang2023}
obtained a time-independent density upper bound and exponential
relaxation in the periodic problem for $\beta>4/3$ and $\gamma>1$.

For the Cauchy problem, the weighted local construction of Li and
Liang~\cite{LiLiang2014} and the global power-law theory of Huang and
Li~\cite{HuangLi2022} provide the starting point. The extension of the
local construction to general coefficients is proved in
Lemma~\ref{ws:lem:general-local}.

The density bound rests on fixed density moments and the dissipation
cost of successive primitive-level crossings. On the torus, pressure
damping and relative-pressure dissipation give uniform moments directly.
In the whole space, a weighted Trudinger estimate first yields
subpolynomial moment growth. A bounded pressure test gives preliminary
energy decay and logarithmically weighted dissipation. These estimates
make the moments uniform in time. Flux estimates and primitive-level
costs then close the density argument in both domains.

\subsection{The periodic problem}\label{per:sec:introduction}

Let $\Omega=\TT$ and normalize the total mass to one.

Set $R_t=\sup_{0\le s\le t}\|\rho(s)\|_{L^{\infty}}$.
For the estimates in both domains, define the effective viscous flux by
\begin{equation}\label{eq:flux-definition}
 G=\begin{cases}
 (2\mu+\lambda(\rho))\diver u-(P-\bar P),&\Omega=\TT,\\
 (2\mu+\lambda(\rho))\diver u-P,&\Omega=\R^2.
 \end{cases}
\end{equation}
Define $A_1,A_2,A_3\ge0$ by
\begin{equation}\label{eq:energies}
\begin{aligned}
 A_1^2(t)&=\int_\Omega\left(\mu\omega^2+
                    \frac{G^2}{2\mu+\lambda(\rho)}\right)\dd,
 &A_2^2(t)&=\int_\Omega\rho|\dot u|^2\dd,\\
 A_3^2(t)&=\int_\Omega\left((2\mu+\lambda(\rho))(\diver u)^2
                       +\mu\omega^2\right)\dd.
\end{aligned}
\end{equation}

Our main results are stated as follows.

\begin{theorem}[Global strong solutions and exponential relaxation]\label{per:thm:main}
Suppose that the structural condition~\ref{per:ass:lambda} holds.
Let $q\in(2,\infty)$ and $C_0\ge1$.
Assume that $0\le\rho_0\in W^{1,q}(\TT)$ and $u_0\in H^1(\TT)$ satisfy
\begin{equation}\label{per:eq:initial}
 \int\rho_0\dd=1,\qquad
 \|\rho_0\|_{L^{\infty}}+\|u_0\|_{H^1}\le C_0.
\end{equation}
Then \eqref{eq:system}, with the initial conditions
\eqref{eq:initial-traces}, admits a unique global strong solution.
For every $T>0$ and $1<s<q/(q-1)$, it has the regularity
\begin{equation}\label{per:eq:strong-class}
\begin{aligned}
 &\rho\in C([0,T];W^{1,q}),\qquad
             \rho_t\in L^\infty(0,T;L^2),\\
 &u\in L^\infty(0,T;H^1)\cap L^2(0,T;H^2)
                         \cap L^s(0,T;W^{2,q}),\\
 &\sqrt t\,u\in L^2(0,T;W^{2,q}),\qquad
             \sqrt t\,u_t\in L^2(0,T;H^1),\\
 &\sqrt\rho\,u_t\in L^2(0,T;L^2),\qquad
             \rho u,\sqrt\rho\,u\in C([0,T];L^2).
\end{aligned}
\end{equation}
The equations hold almost everywhere. The initial density and momentum
are attained in the continuous spaces displayed above.
Moreover, the following bound holds for every $S>0$:
\begin{equation}\label{per:eq:main-bound}
\begin{aligned}
 &\operatorname*{ess\,sup}_{0<t\le S}\bigl(\|\rho(t)\|_{L^{\infty}}+\|u(t)\|_{H^1}^2
                         +A_1^2(t)+A_3^2(t)\bigr)
       \\
 &\qquad+\int_0^S A_2^2\,\mathrm dt
       +\operatorname*{ess\,sup}_{0<t\le S}\min\{1,t\}A_2^2(t)
       +\int_0^S \min\{1,t\}\|\nabla\dot u\|_{L^{2}}^2\,\mathrm dt
 \le C .
\end{aligned}
\end{equation}
The constant is independent of $S$. It depends only on
$\mu,\gamma,C_0$, the constants in the structural
condition~\ref{per:ass:lambda}, and a fixed bound for $\lambda$ on
$[0,2]$. In particular, the low-order bound is independent of
$\|\rho_0\|_{W^{1,q}}$ and of any positive lower density bound used in the
construction. Let
\[
 u_\infty=\int\rho_0u_0\dd.
\]
For every $1\le p<\infty$, there are constants $C_p,\eta_p>0$,
depending only on the quantities listed above and on $p$, such that
\begin{equation}\label{per:eq:main-decay}
 \|\rho(t)-1\|_{L^{p}}+\|u(t)-u_\infty\|_{W^{1,p}}
 \le C_pe^{-\eta_pt},\qquad t\ge1.
\end{equation}
\end{theorem}

\begin{remark}[Comparison with power-law results]\label{per:rem:comparison}
For $\lambda(r)=r^\beta$, Theorem~\ref{per:thm:main} applies to every
$\beta>1$ and $\gamma>1$, including $1<\beta\le4/3$. It extends the
strong existence range of Huang and Li~\cite[Theorem~1.1]{HuangLi2016}.
It also extends the uniform bounds and exponential relaxation of
Fan, Li, and Wang~\cite[Theorem~1.2]{FanLiWang2023}. Both extensions
retain the same initial regularity and allow vacuum.

For $\beta>1$, the weak existence theorem of Li and
Yang~\cite[Theorem~4]{LiYang2026} requires a positive lower bound for the
initial density and an additional higher-integrability bound for the
approximate momenta. For the $W^{1,q}$ density data considered here,
we obtain a strong solution without either requirement.
\end{remark}

\begin{remark}[Higher regularity]\label{per:rem:scope}
The density gradient is controlled on finite time intervals.
Classical regularity requires additional assumptions, since
\eqref{eq:lambda-regularity} does not control $\lambda''$ near vacuum.
Higher coefficient derivatives are used only in the positive-density
construction; the vacuum limit uses local bounds for $\lambda'$.
\end{remark}

\subsection{The Cauchy problem with far-field vacuum}\label{ws:sec:introduction}

Here $\Omega=\R^2$, with vacuum at infinity. We measure velocity in
the homogeneous energy class: an $H^1_{\rm loc}$ representative with
gradient in $L^2$, without an unweighted $L^2$ condition.
Fix $q>2$, $a\in(1,2)$, and any $\eta_0>0$, and set
\begin{equation}\label{ws:eq:weight}
 \bar x=(e+|x|^2)^{1/2}\log^{1+\eta_0}(e+|x|^2).
\end{equation}
The initial data satisfy
\begin{equation}\label{ws:eq:initial-data}
\begin{gathered}
 \rho_0\ge0,\qquad
 0<M:=\int\rho_0\,\dd<\infty,\qquad
 \bar x^a\rho_0\in L^1\cap H^1\cap W^{1,q},\\
 \nabla u_0\in L^2,\qquad \sqrt{\rho_0}u_0\in L^2.
\end{gathered}
\end{equation}
The initial traces prescribed for \eqref{eq:system} are $\rho_0$ and
$\rho_0u_0$. Constants denoted by $C$ in the low-order estimates depend
on $\mu,\gamma,a,\eta_0,M,M^{-1}$, the fixed function $\lambda$, and
\[
 \|\bar x^a\rho_0\|_{L^1}+\|\rho_0\|_{L^\infty}+
 \|\sqrt{\rho_0}u_0\|_{L^2}+\|\nabla u_0\|_{L^2}.
\]
They are independent of the final time and of a positive lower bound
on an approximating density. Dependence on $\lambda$ includes
$C_\lambda$, its $C^1$ bounds on bounded intervals, and the
thresholds in its growth condition at infinity.
Constants $C_T$ may additionally depend
on $T,q$, and $\|\bar x^a\rho_0\|_{H^1\cap W^{1,q}}$.

\begin{theorem}\label{ws:thm:main}
Let $\gamma>1$ and suppose that the structural
condition~\ref{per:ass:lambda} holds.
Data satisfying \eqref{ws:eq:initial-data} generate a unique global strong solution of
\eqref{eq:system}. No compatibility condition involving the initial
acceleration is required. For every finite $T>0$,
\begin{equation}\label{ws:eq:strong-class}
\begin{aligned}
 &\rho\in C([0,T];L^1\cap H^1\cap W^{1,q}),\qquad
       \bar x^a\rho\in L^\infty(0,T;L^1\cap H^1\cap W^{1,q}),\\
 &\sqrt\rho u,\ \nabla u,\ \bar x^{-1}u,\
       \sqrt t\,\sqrt\rho\,u_t\in L^\infty(0,T;L^2),\\
 &\nabla u\in L^2(0,T;H^1)
                \cap L^{(q+1)/q}(0,T;W^{1,q}),\qquad
       \sqrt t\,\nabla u\in L^2(0,T;W^{1,q}),\\
 &\sqrt\rho\,u_t,\ \sqrt t\,\nabla u_t,\
       \sqrt t\,\bar x^{-1}u_t\in L^2((0,T)\times\R^2),\\
 &\rho u,\ \sqrt\rho\,u\in C([0,T];L^2).
\end{aligned}
\end{equation}
Moreover, $\nabla u\in L^1(0,T;L^\infty)$, and both mass and total
momentum are conserved. The following uniform bound holds:
\begin{equation}\label{ws:eq:uniform-core}
\begin{aligned}
 &\sup_{t\ge0}\bigl(\|\rho(t)\|_{L^\infty}+
                   \|\nabla u(t)\|_{L^2}^2+A_1^2(t)\bigr)\\
 &\quad+\int_0^\infty
       \bigl(A_3^2+\|\sqrt\rho\,\dot u\|_{L^2}^2+\|P\|_{L^2}^2\bigr)\,dt\\
 &\quad+\sup_{t>0}\min\{1,t\}
       \|\sqrt{\rho(t)}\,\dot u(t)\|_{L^2}^2
       +\int_0^\infty\min\{1,t\}\|\nabla\dot u(t)\|_{L^2}^2\,dt\le C.
\end{aligned}
\end{equation}
For every $\tau>0$,
\begin{equation}\label{ws:eq:flux-square-global}
 \int_\tau^\infty\|G(t)\|_{L^\infty}^2\,dt\le C_\tau.
\end{equation}
Define the conserved mean velocity by
\[
 u_\infty=\frac1M\int\rho_0u_0\dd .
\]
For $t\ge1$, the following decay estimates hold:
\begin{equation}\label{ws:eq:decay-density-pressure}
\begin{aligned}
 \|\rho(t)\|_{L^p}
 &\le C_p(1+t)^{-\frac1\gamma(1-\frac1p)},
       &&1<p\le\infty,\\
 \|P(t)\|_{L^r}
 &\le C_r(1+t)^{-1+\frac1{\gamma r}},
       &&1\le r\le\infty.
\end{aligned}
\end{equation}
Moreover,
\begin{equation}\label{ws:eq:decay-velocity}
\begin{aligned}
 \|\nabla u(t)\|_{L^p}
 &\le C_p(1+t)^{-1+\frac1{\gamma p}},
       &&2\le p<\infty,\\
 \|\sqrt{\rho(t)}\,\dot u(t)\|_{L^2}
 &\le C(1+t)^{-2+\frac1{2\gamma}}.
\end{aligned}
\end{equation}
For the same conserved mean velocity $u_\infty$,
\begin{equation}\label{ws:eq:decay-kinetic}
\begin{aligned}
 \int\rho|u-u_\infty|^2\dd
 &\le C(1+t)^{-2+\frac1\gamma},\\
 \frac12\int\rho|u-u_\infty|^2\dd
       +\frac1{\gamma-1}\int P\dd
 &\le C(1+t)^{-1+\frac1\gamma}.
\end{aligned}
\end{equation}
The constants are independent of the terminal time. The matching
lower bounds and the precise sharpness statements are given in
Proposition~\ref{ws:prop:decay-sharpness}.
\end{theorem}

\begin{remark}[Coefficient examples and scope]\label{rem:joint-scope}
For $\lambda_0\ge0$, the following families satisfy the common
structural condition in both domains:
\[
\begin{array}{c|c}
 \lambda(r)-\lambda_0&\text{admissible range}\\ \hline
 r^\beta[\log(e+r)]^s&\beta>1,\ s\in\R\\
 r[\log(e+r)]^\alpha&\alpha>1/2\\
 r\log(e+r)[\log\log(e^e+r)]^\delta&\delta\in\R
\end{array}
\]
Positive finite sums are admissible. Multiplication by a positive
bounded factor bounded away from zero also preserves admissibility,
provided its logarithmic derivative is bounded and the stated
regularity is preserved. For example, $2+\sin(4\log(e+r))$ gives
nonmonotone coefficients.
Indeed, $\lambda(r)\ge cr[\log(e+r)]^\alpha$ for large $r$ implies
\[
 \theta(r)\ge\int_{r/2}^r\frac{\lambda(z)}z\,\mathrm dz
       \ge c'r[\log(e+r)]^\alpha.
\]
Thus \eqref{per:eq:primitive-growth} holds for $\alpha>1/2$, with
$\kappa=\min\{\alpha,1\}$. The range $\kappa\le1/2$ is not covered.
\end{remark}

\begin{remark}[The common growth range]\label{rem:domain-growth}
The periodic and Cauchy results use the same viscosity growth condition,
including $\lambda(r)\asymp r[\log(e+r)]^\alpha$ for every $\alpha>1/2$.
In the whole space, preliminary energy decay is obtained before the
density upper bound; see Lemma~\ref{ws:lem:preliminary-energy}.
This decay makes the coefficient in \eqref{ws:eq:log-time-moment-ode}
integrable over the whole lifespan. It yields uniform density moments
and pressure dissipation without requiring a fixed ball to contain a
prescribed mass for all time. The flux bound and primitive-level costs
then apply as in the periodic argument. The different equilibria lead
to exponential relaxation on the torus and algebraic relaxation toward
vacuum in the whole space.
\end{remark}

Sections~\ref{sec:periodic-density} and~\ref{sec:cauchy-density}
give the density bounds. Construction and uniqueness are proved in
Section~\ref{sec:strong-solutions}, and large-time behavior in
Section~\ref{sec:large-time}.

Constants $C$ may change from line to line. In low-order estimates,
their dependence is stated in the corresponding main theorem; they are
independent of the terminal time and positive density lower bounds.
Constants $C_T$ may depend on the fixed finite time. In strong and
uniqueness estimates, they may also depend on higher initial norms.
The constants $C_\tau$ may depend on a fixed positive time $\tau$.
Other subscripts indicate additional parameters.

\section{Uniform density estimates in the periodic problem}\label{sec:periodic-density}

Throughout this section, $\Omega=\TT$ and Assumption~\ref{per:ass:lambda} is imposed.

\subsection{Uniform density moments and pressure dissipation}\label{per:sec:moment}

\begin{proposition}[Uniform low-order estimates]\label{per:prop:apriori}
Let $(\rho,u)$ be a smooth positive solution of \eqref{eq:system}
on $[0,S]$ satisfying \eqref{per:eq:initial}. Then \eqref{per:eq:main-bound}
holds with the dependence stated in Theorem~\ref{per:thm:main}.
\end{proposition}
All estimates below concern these solutions; the proof is completed in
Subsection~\ref{per:sec:density}.

Conservation of momentum allows the Galilean change
\[
 (\rho,u)(t,x)\longmapsto
 \bigl(\rho(t,x+tu_\infty),\,u(t,x+tu_\infty)-u_\infty\bigr).
\]
It preserves the equations, density norms, and $A_1,A_2,A_3$.
Moreover, $|u_\infty|\le C(C_0)$. We henceforth use the zero-momentum
frame $\int\rho u\dd=0$; see \cite[Appendix~I]{FanLiWang2023}.

Testing \eqref{eq:system}$_2$ by $u$ and using
\eqref{eq:system}$_1$ and periodicity gives
\[
 \frac{\mathrm d}{\mathrm dt}\int
 \left(\frac12\rho|u|^2+\frac{\rho^\gamma}{\gamma-1}\right)\dd
 +A_3^2=0.
\]
Here $\int|\nabla u|^2=\int((\diver u)^2+\omega^2)$. Also,
\[
 |\bar u|=\left|\int\rho(u-\bar u)\dd\right|
 \le\|\rho\|_{L^{\gamma}}\|u-\bar u\|_{L^{\gamma/(\gamma-1)}}
 \le C\|\nabla u\|_{L^{2}}.
\]
Thus
\begin{equation}\label{per:eq:basic}
 \sup_{0\le t\le S}\int(\rho|u|^2+\rho^\gamma)\dd
       +\int_0^S A_3^2\,\mathrm dt\le C,
 \qquad \|u\|_{H^1}\le CA_3.
\end{equation}

\begin{lemma}[Coefficient bounds]\label{lem:coefficient}
Under \eqref{eq:lambda-regularity}, let $b=\max\{1,C_\lambda\}$. Then
\begin{equation}\label{eq:coefficient-bounds}
\begin{gathered}
 2\mu+\lambda(r)\le C(1+r^b),\qquad
 \lambda(r)\le C(1+\theta(r)_+),\qquad r\ge0,\\
 2^{-b}(2\mu+\lambda(r))\le2\mu+\lambda(z)
       \le2^b(2\mu+\lambda(r)),\qquad r/2\le z\le r.
\end{gathered}
\end{equation}
Here $\theta(0)_+=0$. These bounds apply in both domains.
\end{lemma}
\begin{proof}
Integrate $|(\log(2\mu+\lambda(r)))'|\le C_\lambda/r$.
For $r\ge2$, the dyadic comparison also gives
\[
 \theta(r)\ge\int_{r/2}^r\frac{2\mu+\lambda(z)}z\,\mathrm dz
       \ge2^{-b}\log2\,(2\mu+\lambda(r)).
\]
The coefficient bound on $[0,2]$ treats the remaining interval.
\end{proof}

Let $k$ be an integer satisfying
\begin{equation}\label{per:eq:k-choice}
 k\ge3,\qquad (2\kappa-1)k\ge1.
\end{equation}
Constants may depend on this fixed $k$.

\begin{lemma}[Uniform density moments and tails]\label{per:lem:moment}
For every $k$ satisfying \eqref{per:eq:k-choice},
\begin{equation}\label{per:eq:moment}
 \sup_{0\le t\le S}\int\left(
 \rho\theta(\rho)_+^k+
 \rho^{k+1}[\log(e+\rho)]^{\kappa k}\right)\dd\le C_k.
\end{equation}
Consequently, for every fixed $1\le p<\infty$,
\begin{equation}\label{per:eq:uniform-density-moments}
 \sup_{0\le t\le S}\|\rho(t)\|_{L^{p}}\le C_p.
\end{equation}
For every fixed $p\ge2$, there are $L_p\ge2$ and $C_p$ such that
\begin{equation}\label{per:eq:density-tail-integrability}
 \int_0^S\int(\rho-L_p)_+^p\dd\,\mathrm dt\le C_p.
\end{equation}
All these constants are independent of $S$.
\end{lemma}
\begin{proof}
Put $\psi=(-\Delta)^{-1}\diver(\rho u)$ with zero mean.
For $R_i=(-\Delta)^{-1/2}\partial_i$ and
$[u_i,R_iR_j]f=u_iR_iR_jf-R_iR_j(u_if)$, define
\[
 F=\sum_{i,j=1}^2[u_i,R_iR_j](\rho u_j).
\]
For smooth positive solutions, the continuity equation gives
\begin{equation}\label{eq:theta-flow}
 \Dt\theta(\rho)+P=
 \begin{cases}\bar P-G,&\Omega=\TT,\\-G,&\Omega=\R^2.
 \end{cases}
\end{equation}
Its use at vacuum is justified through positive approximations.

Combining \eqref{eq:theta-flow} with the conservative momentum equation gives
\begin{equation}\label{per:eq:transport-potential}
 \Dt(\theta(\rho)-\psi)+P=\bar P-\bar G-F.
\end{equation}
Here $G-\bar G+\Dt\psi=F$, and all inverse Laplacians act on zero-mean functions.

For $M>0$ to be fixed, set
\[
 \varphi=(\theta(\rho)-\psi-M)_+,\qquad
 \mathcal M_k(t)=\int\rho\varphi^k\dd.
\]
The growth assumption gives
\[
 \int\rho^{k+1}[\log(e+\rho)]^{\kappa k}\dd
 \le C_k\left(1+\int\rho\theta(\rho)_+^k\dd\right).
\]
The following potential estimates hold in either domain whenever mass and
kinetic energy are bounded. Hardy--Littlewood--Sobolev and mass interpolation give
\begin{equation}\label{eq:potential-comparison}
\begin{aligned}
 \|\psi\|_{L^{2(k+1)}}
 &\le C_k\|\rho u\|_{L^{2(k+1)/(k+2)}}
       \le C_k\|\rho\|_{L^{k+1}}^{1/2},\\
 \int\rho|\psi|^k\dd
 &\le\|\rho\|_{L^{2(k+1)/(k+2)}}
                           \|\psi\|_{L^{2(k+1)}}^k\\
 &\le C_k\|\rho\|_{L^{k+1}}^{(k+1)/2}
       \le C_k\left(1+\int\rho\theta(\rho)_+^k\dd\right)^{1/2}.
\end{aligned}
\end{equation}
Since $\theta(\rho)_+\le\varphi+|\psi|+M$, absorbing the square root yields
\begin{equation}\label{per:eq:moment-comparison}
 \int\left(\rho\theta(\rho)_+^k+
           \rho^{k+1}[\log(e+\rho)]^{\kappa k}\right)\dd
 \le C_{k,M}(1+\mathcal M_k).
\end{equation}
In particular,
\begin{equation}\label{per:eq:potential-mass-moment}
 \int\rho|\psi|^k\dd
 \le C_{k,M}(1+\mathcal M_k)^{1/2}.
\end{equation}

For the momentum estimate, \eqref{per:eq:basic} and Trudinger's
inequality~\cite{Trudinger1967} give
\[
 \int\exp\left(\frac{c|u|^2}{A_3^2}\right)\dd\le C
 \quad\text{when }A_3>0;
 \qquad u=0\quad\text{when }A_3=0.
\]
Orlicz--Young then gives
\[
 \int\rho^{k+1}|u|^{k+1}\dd
 \le C_kA_3^{k+1}\left(
 1+\int\rho^{k+1}[\log(e+\rho)]^{(k+1)/2}\dd\right).
\]
Since $\kappa k\ge(k+1)/2$, \eqref{per:eq:moment-comparison} implies
\begin{equation}\label{per:eq:momentum-fixed}
 \|\rho u\|_{L^{k+1}}
 \le C_{k,M}A_3(1+\mathcal M_k)^{1/(k+1)}.
\end{equation}
H\"older and the commutator theorem~\cite{CRW1976} give
\begin{equation}\label{eq:moment-commutator}
\begin{aligned}
 \int\rho|F|\varphi^{k-1}\dd
 &\le\mathcal M_k^{1-1/k}\|\rho\|_{L^{k+1}}^{1/k}\|F\|_{L^{k+1}}\\
 &\le C_{k,M}\|\nabla u\|_{L^2}\|\rho u\|_{L^{k+1}}
                         (1+\mathcal M_k)^{k/(k+1)}.
\end{aligned}
\end{equation}
The estimate also holds for the whole-space commutator with zero
truncation shift. On the torus, \eqref{per:eq:momentum-fixed} bounds the
commutator term by $C_{k,M}A_3^2(1+\mathcal M_k)$.

For the mean, \eqref{eq:coefficient-bounds} gives
\[
\begin{aligned}
 |\bar G|^2
 &=\left|\int\lambda\diver u\dd\right|^2
 \le A_3^2\int\lambda(\rho)\dd\\
 &\le C_kA_3^2\left(
       1+\left[\int\rho\theta(\rho)_+^k\dd\right]^{1/k}\right)
 \le C_{k,M}A_3^2(1+\mathcal M_k)^{1/k}.
\end{aligned}
\]
Using $|\bar G|\le|\bar G|^2+1$ and
$\int\rho\varphi^{k-1}\dd\le\mathcal M_k^{1-1/k}$, we obtain
\[
 |\bar G|\int\rho\varphi^{k-1}\dd
 \le C_{k,M}A_3^2(1+\mathcal M_k)
       +\int\rho\varphi^{k-1}\dd.
\]
Thus testing \eqref{per:eq:transport-potential} gives
\[
 \mathcal M_k'+k\int\rho^{\gamma+1}\varphi^{k-1}\dd
 \le C_{k,M}A_3^2(1+\mathcal M_k)
                    +Ck\int\rho\varphi^{k-1}\dd.
\]
Here $C$ is independent of $M$, since $\bar P\le C$ and the mean
estimate contributes coefficient one. Choose
\[
 M>1+\max\left\{\theta\bigl((4C)^{1/\gamma}\bigr),0\right\}.
\]
On $\{\rho^\gamma\ge4C\}$, the last term is absorbed by half
the pressure damping. On $\{\rho^\gamma<4C,\varphi>0\}$,
we have $\varphi\le|\psi|$. The potential estimate gives
\[
 \|\psi\|_{L^{\infty}}
 \le C\|\rho u\|_{L^{8/3}}
 \le C\|\rho\|_{L^{4}}\|u\|_{L^{8}}
 \le C_{k,M}A_3(1+\mathcal M_k)^{1/(k+1)}.
\]
By \eqref{per:eq:potential-mass-moment} and $k\ge3$,
\[
\begin{aligned}
 \int_{\{\rho^\gamma<4C\}}\rho\varphi^{k-1}\dd
 &\le\|\psi\|_{L^{\infty}}^2\int\rho|\psi|^{k-3}\dd\\
 &\le C_{k,M}A_3^2(1+\mathcal M_k)^{2/(k+1)+1/2}
 \le C_{k,M}A_3^2(1+\mathcal M_k).
\end{aligned}
\]
We conclude that
\begin{equation}\label{per:eq:uniform-moment-differential}
 \mathcal M_k'+\frac k2\int\rho^{\gamma+1}\varphi^{k-1}\dd
 \le C_{k,M}A_3^2(1+\mathcal M_k).
\end{equation}
Since $\mathcal M_k(0)\le C_k$, Gronwall and \eqref{per:eq:basic} yield
\[
 \sup_{0\le t\le S}\mathcal M_k(t)
 +\int_0^S\int\rho^{\gamma+1}\varphi^{k-1}\dd\,\mathrm dt
 \le C_k.
\]
By \eqref{per:eq:moment-comparison}, this proves \eqref{per:eq:moment}
and \eqref{per:eq:uniform-density-moments}.

For the tail estimate, fix $p\ge3$ and choose $k\ge p$ in the
preceding argument. Choose $L_p\ge2$ so that
$\theta(r)\ge2M$ and $\theta(r)\ge cr$ for $r\ge L_p$.
On $\{\rho\ge L_p\}$, either
$\varphi\ge\theta(\rho)/4$, in which case
$\rho^p\le C\rho^{\gamma+1}\varphi^{k-1}$, or
$\psi\ge\theta(\rho)/4\ge c\rho$. Consequently,
\[
\begin{aligned}
 \int_{\{\rho\ge L_p\}}\rho^p\dd
 &\le C_p\int\rho^{\gamma+1}\varphi^{k-1}\dd
       +C_p\|\psi\|_{L^{\infty}}^2\int\rho|\psi|^{p-3}\dd\\
 &\le C_p\int\rho^{\gamma+1}\varphi^{k-1}\dd+C_pA_3^2.
\end{aligned}
\]
The last step uses \eqref{per:eq:potential-mass-moment} and the bound
for $\mathcal M_k$. Integrating proves the tail estimate for $p\ge3$.
For $2\le p<3$, take $L_p=L_3+1$ and use
$(\rho-L_3-1)_+^p\le(\rho-L_3)_+^3$.
\end{proof}

\begin{lemma}[Global relative-pressure dissipation]\label{per:lem:pressure-dissipation}
There is a constant independent of $S$ such that
\begin{equation}\label{per:eq:pressure-dissipation}
\begin{aligned}
 &\int_0^S\int(\rho+1)^{\gamma-1}(\rho-1)^2\dd\,\mathrm dt\\
 &\qquad+\int_0^S
       \bigl(\|P-1\|_{L^{2}}^2+\|P-1\|_{L^{4}}^4\bigr)\,\mathrm dt\le C.
\end{aligned}
\end{equation}
\end{lemma}
\begin{proof}
Test momentum by $\nabla\Delta^{-1}(\rho-1)$, as in
\cite[Propositions~3.3 and~5.4]{FanLiWang2023}. Continuity gives
\begin{equation}\label{per:eq:pressure-test}
\begin{aligned}
 \int(P-1)(\rho-1)\dd
 ={}&\frac{\mathrm d}{\mathrm dt}
       \int\rho u\cdot\nabla\Delta^{-1}(\rho-1)\dd
       +\|\nabla\Delta^{-1}\diver(\rho u)\|_{L^{2}}^2\\
 &-\int\rho u\otimes u:\nabla^2\Delta^{-1}(\rho-1)\dd
       +\int(2\mu+\lambda)\diver u\,(\rho-1)\dd.
\end{aligned}
\end{equation}
The uniform moments, elliptic estimates, and \eqref{per:eq:basic} imply
\[
\begin{aligned}
 \|\nabla\Delta^{-1}\diver(\rho u)\|_{L^{2}}^2
 &\le C\|\rho\|_{L^{4}}^2\|u\|_{L^{4}}^2\le CA_3^2,\\
 \left|\int\rho u\otimes u:
                    \nabla^2\Delta^{-1}(\rho-1)\dd\right|
 &\le C\|\rho\|_{L^{4}}\|u\|_{L^{4}}^2\|\rho-1\|_{L^{4}}\le CA_3^2,\\
 \left|\int\rho u\cdot\nabla\Delta^{-1}(\rho-1)\dd\right|
 &\le\|\sqrt\rho\,u\|_{L^{2}}\|\rho\|_{L^{2}}^{1/2}
              \|\nabla\Delta^{-1}(\rho-1)\|_{L^{4}}\le C.
\end{aligned}
\]
Splitting at $\rho=2L_{b+2}$ and using \eqref{eq:coefficient-bounds} gives
\[
\begin{aligned}
 \left|\int(2\mu+\lambda)\diver u\,(\rho-1)\dd\right|
 &\le\varepsilon\|\rho-1\|_{L^{2}}^2+C_\varepsilon A_3^2
       +C\int_{\{\rho>2L_{b+2}\}}\rho^{b+2}\dd\\
 &\le\varepsilon\|\rho-1\|_{L^{2}}^2+C_\varepsilon A_3^2
       +C\int(\rho-L_{b+2})_+^{b+2}\dd.
\end{aligned}
\]
Since
\[
 (\rho^\gamma-1)(\rho-1)
 \ge c_\gamma(\rho+1)^{\gamma-1}(\rho-1)^2,
\]
integration of \eqref{per:eq:pressure-test}, with $\varepsilon<c_\gamma/2$,
and \eqref{per:eq:density-tail-integrability} prove the first assertion.
For $p=2,4$,
\[
 |\rho^\gamma-1|^p
 \le C_p(\rho-1)^2+C_p(\rho-L_{p\gamma})_+^{p\gamma}.
\]
The first assertion and \eqref{per:eq:density-tail-integrability}
therefore prove the remaining integrals.
\end{proof}

The definition of $G$ and the preceding estimates yield
\begin{equation}\label{per:eq:energy-comparison}
\begin{aligned}
 A_1^2&\le C\bigl(A_3^2+\|P-\bar P\|_{L^{2}}^2\bigr),\qquad
 A_3^2\le C\bigl(A_1^2+\|P-\bar P\|_{L^{2}}^2\bigr),\\
 e+A_1^2&\asymp1+A_3^2,\qquad
 \int_0^S A_1^2\,\mathrm dt\le C.
\end{aligned}
\end{equation}
These constants are independent of $S$ and $R_S$.

\subsection{A weighted effective-flux estimate}\label{per:sec:flux}

\begin{lemma}[Flux estimates in both domains]\label{lem:flux-elliptic}
For $\Omega=\TT$ or $\R^2$, $1<p<\infty$, and
$R\ge\max\{1,\|\rho\|_{L^\infty}\}$,
\begin{equation}\label{eq:flux-elliptic}
\begin{gathered}
 \rho\dot u=\nabla G+\mu\nabla^\perp\omega,\qquad
 \|\nabla G\|_{L^p}+\|\nabla\omega\|_{L^p}
       \le C_p\|\rho\dot u\|_{L^p},\\
 \|\nabla G\|_{L^2}+\|\nabla\omega\|_{L^2}
       \le C\sqrt R\,A_2,\qquad
 \|G\|_{L^2}\le CR^{b/2}A_1.
\end{gathered}
\end{equation}
The constants are independent of $R$ and time.
\end{lemma}
\begin{proof}
Take divergence and curl of momentum and use Calder\'on--Zygmund estimates.
The $L^2$ bounds follow from $\|\rho\dot u\|_{L^2}\le\sqrt R\,A_2$
and \eqref{eq:coefficient-bounds}.
\end{proof}

\begin{lemma}[Logarithmic energy]\label{lem:log-energy}
Let $(\rho,u)$ be a local strong solution on $[0,S]$, with
$\Omega=\TT$ or $\R^2$ and $\lambda$ satisfying \eqref{eq:lambda-regularity}.
Write $\bar P=0$ when $\Omega=\R^2$. Suppose that, for some $K\ge1$,
the following bounds hold on $[0,S]$:
\[
\begin{gathered}
 \sup_{t\le S}\|P\|_{L^4}
 +\mathop{\mathrm{ess\,limsup}}_{t\downarrow0}A_1^2(t)
 +\int_0^S\bigl(A_3^2+\|P-\bar P\|_{L^2}^2\bigr)\,dt\le K,\\
 A_3^2\le K(e+A_1^2).
\end{gathered}
\]
On the torus, also assume $|\bar G|\le K A_3$.
For $R\ge\max\{1,\sup_{t\le S}\|\rho\|_{L^\infty}\}$,
\begin{equation}\label{eq:log-energy}
 \mathop{\mathrm{ess\,sup}}_{0<t\le S}\log(e+A_1^2(t))
 +\int_0^S\frac{A_2^2}{e+A_1^2}\,dt
 \le C_K R\log(e+R).
\end{equation}
The constant depends on $K$ and the fixed coefficients, and is
independent of $S$, $R$, and a positive density lower bound.
\end{lemma}
\begin{proof}
Test momentum by $2\dot u$ and use
\[
\begin{aligned}
 \diver\dot u&=\Dt\diver u+(\diver u)^2-2\det\nabla u,\\
 \curl\dot u&=\Dt\omega+\omega\diver u,\qquad
 \Dt\lambda=-\rho\lambda'\diver u,\qquad
 \Dt P=-\gamma P\diver u.
\end{aligned}
\]
The coefficient bound, \eqref{eq:flux-elliptic}, vorticity interpolation,
and the Hardy--BMO estimate for the Jacobian \cite{CLMS1993} give
\[
\begin{aligned}
 (A_1^2)'+2A_2^2
 \le{}&C A_3\left(
       \left\|\frac{G^2}{2\mu+\lambda}\right\|_{L^2}
       +\left\|\frac{G(P+\bar P)}{2\mu+\lambda}\right\|_{L^2}\right)\\
 &+C\sqrt R\,A_2(A_1A_3+A_3^2)
             +C_K\bigl(A_3^2+\|P-\bar P\|_{L^2}^2\bigr).
\end{aligned}
\]
Here the last term includes the extra mean-pressure term on the torus:
\[
 \bar P_t=-(\gamma-1)\int(P-\bar P)\diver u\dd,\qquad
 \int\frac{G}{2\mu+\lambda}\dd
       =-\int\frac{P-\bar P}{2\mu+\lambda}\dd,
\]
so its absolute value is bounded by
$C_K A_3\|P-\bar P\|_{L^2}^2
\le C_K(A_3^2+\|P-\bar P\|_{L^2}^2)$.
This term is absent in the whole space; no whole-space mean of $G$ is used.

For $0<\delta\le1/4$, we have in both domains
\begin{equation}\label{eq:weighted-quadratic}
 \left\|\frac{G^2}{2\mu+\lambda}\right\|_{L^2}
 \le\frac C{\sqrt\delta}R^{(1+b\delta)/2}A_1A_2
              +C_KA_1^{1-\delta}A_3^{1+\delta}.
\end{equation}
Indeed, H\"older's inequality bounds the left-hand side by
$CA_1^{1-\delta}\|G\|_{L^{2(1+\delta)/\delta}}^{1+\delta}$.
Quantitative Gagliardo--Nirenberg gives
\[
 \|G-\bar G\|_{L^{2(1+\delta)/\delta}}^{1+\delta}
 \le\frac C{\sqrt\delta}
         \|G-\bar G\|_{L^2}^{\delta}\|\nabla G\|_{L^2}
 \qquad(\Omega=\TT),
\]
and the same inequality with $G$ on $\R^2$.
Use \eqref{eq:flux-elliptic} and, on the torus, $|\bar G|\le K A_3$.
The interpolation constant is uniform since
$\delta^{-(1+\delta)/2}\le C\delta^{-1/2}$.
The last term in \eqref{eq:weighted-quadratic} is unnecessary on $\R^2$,
but the displayed bound remains valid there.

The pressure bound and Young's inequality imply
\[
\begin{aligned}
 A_3\left\|\frac{G^2}{2\mu+\lambda}\right\|_{L^2}
 &\le\varepsilon A_2^2+
     \frac{C_{K,\varepsilon}}\delta R^{1+b\delta}
                         (e+A_1^2)A_3^2,\\
 A_3\left\|\frac{G(P+\bar P)}{2\mu+\lambda}\right\|_{L^2}
 &\le C_K A_3
        \left\|\frac{G^2}{2\mu+\lambda}\right\|_{L^2}^{1/2}\\
 &\le\varepsilon A_2^2+
     \frac{C_{K,\varepsilon}}\delta R^{1+b\delta}(A_1^2+A_3^2).
\end{aligned}
\]
Absorb the acceleration terms, divide by $e+A_1^2$, and use the elementary
bound $A_1^2\le C(A_3^2+\|P-\bar P\|_{L^2}^2)$ to obtain
\begin{equation}\label{eq:log-energy-differential}
 \bigl(\log(e+A_1^2)\bigr)'
       +c\frac{A_2^2}{e+A_1^2}
 \le\frac{C_K}\delta R^{1+b\delta}
                   \bigl(A_3^2+\|P-\bar P\|_{L^2}^2\bigr).
\end{equation}
Integrate from admissible positive times, use the initial upper trace,
and take $\delta=[4\log(e+R)]^{-1}$ to obtain \eqref{eq:log-energy}.
\end{proof}

On the torus, the hypotheses follow from Lemmas~\ref{per:lem:moment}
and \ref{per:lem:pressure-dissipation} and
\eqref{per:eq:energy-comparison}, with $K$ depending only on the
low-order data. We use \eqref{eq:log-energy} with $R=R_S$ below.

\begin{remark}[Quantitative improvement of the logarithmic estimate]
\label{per:rem:log-energy-comparison}
For $\lambda(\rho)=\rho^\beta$, the finite-time logarithmic estimates
of Huang and Li~\cite[(3.52)]{HuangLi2016} and Li and
Yang~\cite[(65)]{LiYang2026} have right-hand sides
$C(\alpha,T)R_S^{1+\alpha\beta}$ for fixed $\alpha\in(0,1)$;
the dependence as $\alpha\downarrow0$ is not quantified there.
Here the explicit factor $\delta^{-1}$ gives the logarithmic correction,
while the uniform moments and pressure dissipation remove final-time
dependence; cf.~\cite[Section~3.3 and Appendix~I]{FanLiWang2023}.
The resulting flux bound in Proposition~\ref{per:prop:flux} gives the
level costs for the near-linear coefficients in Remark~\ref{rem:joint-scope}.
\end{remark}

For either domain, commuting momentum with the material derivative gives
\begin{equation}\label{eq:material-momentum}
\begin{aligned}
 \rho\Dt\dot u_i={}&\mu\Delta\dot u_i
                   +\partial_i\bigl((\mu+\lambda)\diver\dot u\bigr)\\
 &+\mu\partial_j\bigl((\diver u)\partial_j u_i
       -\partial_j u_\ell\partial_\ell u_i
       -\partial_\ell u_j\partial_\ell u_i\bigr)\\
 &+\partial_i\bigl(2(\mu+\lambda)\det\nabla u
       -\rho\lambda'(\diver u)^2+(\gamma-1)P\diver u\bigr)\\
 &-\partial_j\bigl(((\mu+\lambda)\diver u-P)\partial_i u_j\bigr).
\end{aligned}
\end{equation}
Testing by $2\dot u$ and using
$|\rho\lambda'|\le C(\mu+\lambda)$ gives
\begin{equation}\label{eq:material-acceleration}
\begin{aligned}
 &(A_2^2)'+c\|\nabla\dot u\|_{L^2}^2
             +c\int(\mu+\lambda)(\diver\dot u)^2\dd\\
 &\quad\le C\int\bigl[(1+\lambda)|\nabla u|^4
       +P^2|\nabla u|^2
       +\bigl((\mu+\lambda)\diver u-P\bigr)^2|\nabla u|^2\bigr]\dd.
\end{aligned}
\end{equation}
Scalar errors are absorbed against
$\int(\mu+\lambda)(\diver\dot u)^2$, and tensor errors against
$\|\nabla\dot u\|_{L^2}^2$. The constants depend only on
$\mu,\gamma,C_\lambda$.

\begin{lemma}[Uniform positive-time acceleration bound]\label{per:lem:acceleration}
For every $0<\tau<S$,
\begin{equation}\label{per:eq:log-acceleration}
 \sup_{\tau\le t\le S}\log(e+A_2^2(t))
 +\int_\tau^S\frac{\|\nabla\dot u\|_{L^{2}}^2}{e+A_2^2}\,\mathrm dt
 \le C_\tau R_S^{b+1}\log(e+R_S),
\end{equation}
where $C_\tau$ is independent of $S$.
\end{lemma}
\begin{proof}
On the torus,
$(\mu+\lambda)\diver u-P=G-\bar P-\mu\diver u$ and $\bar P\le C$.
H\"older's and Young's inequalities therefore bound the right-hand side of
\eqref{eq:material-acceleration} by
\[
 C\bigl(R_S^b\|\nabla u\|_{L^4}^4+\|G\|_{L^4}^4
                    +A_3^2+\|P-\bar P\|_{L^4}^4\bigr).
\]

To estimate the fourth powers, first use
\[
\begin{aligned}
 \bar G\int\frac1{2\mu+\lambda}\dd
 &=-\int\frac{G-\bar G+P-\bar P}{2\mu+\lambda}\dd,\\
 \int\frac1{2\mu+\lambda}\dd
 &\ge\left(\int(2\mu+\lambda)\dd\right)^{-1}\ge c.
\end{aligned}
\]
Together with \eqref{eq:flux-elliptic} and $|\bar G|\le CA_3$, this yields
\[
\begin{aligned}
 |\bar G|&\le C\bigl(R_S^{1/2}A_2+\|P-\bar P\|_{L^{2}}\bigr),\\
 A_1^2&\le CR_S A_2^2+C\|P-\bar P\|_{L^{2}}^2,\\
 |\bar G|^4&\le CA_3^2
                    \bigl(R_S A_2^2+\|P-\bar P\|_{L^{2}}^2\bigr).
\end{aligned}
\]
Hence \eqref{per:eq:energy-comparison} gives
\[
 A_3^2+\|P-\bar P\|_{L^{2}}^2\le CR_S(e+A_2^2).
\]
Use \eqref{eq:weighted-quadratic} with
$\delta=[4\log(e+R_S)]^{-1}$. Hodge and vorticity interpolation give
\[
 \|\nabla u\|_{L^{4}}^2
 \le C\sqrt{R_S\log(e+R_S)}\,A_1A_2
       +C\bigl(A_3^2+\|P-\bar P\|_{L^{2}}^2+\|P-\bar P\|_{L^{4}}^2\bigr).
\]
Squaring, and applying zero-mean interpolation to $G-\bar G$, yields
\begin{equation}\label{per:eq:velocity-fourth-powers}
\begin{aligned}
 \|\nabla u\|_{L^{4}}^4
 &\le CR_S\log(e+R_S)
        \bigl(A_3^2+\|P-\bar P\|_{L^{2}}^2\bigr)(e+A_2^2)
                             +C\|P-\bar P\|_{L^{4}}^4,\\
 \|G\|_{L^{4}}^4
 &\le CR_S^{b+1}
        \bigl(A_3^2+\|P-\bar P\|_{L^{2}}^2\bigr)(e+A_2^2).
\end{aligned}
\end{equation}
Substitution gives
\begin{equation}\label{per:eq:acceleration-differential}
\begin{aligned}
 (A_2^2)'+c\|\nabla\dot u\|_{L^{2}}^2
 \le{}&CR_S^{b+1}\log(e+R_S)(e+A_2^2)\\
 &\times\bigl(A_3^2+\|P-\bar P\|_{L^{2}}^2+\|P-\bar P\|_{L^{4}}^4\bigr).
\end{aligned}
\end{equation}
By \eqref{per:eq:basic} and \eqref{per:eq:pressure-dissipation}, the last factor
is integrable independently of $S$. Choose $s_\tau\in[\tau/2,\tau]$ with
\[
 A_2^2(s_\tau)\le\frac2\tau\int_0^\tau A_2^2\,\mathrm dt,
 \qquad
 \log(e+A_2^2(s_\tau))\le C_\tau R_S\log(e+R_S),
\]
using \eqref{eq:log-energy} for the second bound. Divide
\eqref{per:eq:acceleration-differential} by $e+A_2^2$ and integrate.
\end{proof}

\begin{lemma}[Periodic interpolation]\label{per:lem:interpolation}
Let $\bar f=0$ and $p,m\ge8$, and set
\[
 \sigma=\frac{2m}{p(m-2)+2m},\qquad
 \alpha=\left(1-\frac2m\right)\sigma.
\]
Then
\begin{equation}\label{per:eq:interpolation}
 \|f\|_{L^{\infty}}
 \le C\sqrt p\,\|f\|_{L^{2}}^\alpha
       \|\nabla f\|_{L^{2}}^{1-\alpha-\sigma}\|\nabla f\|_{L^{m}}^\sigma,
\end{equation}
with $C$ independent of $p,m$.
\end{lemma}
\begin{proof}
Ball averaging and Morrey's inequality give, uniformly for $p,m\ge8$,
\[
 |f(x)|\le C\ell^{-2/p}\|f\|_{L^{p}}
                 +C\ell^{1-2/m}\|\nabla f\|_{L^{m}},
 \qquad 0<\ell\le\ell_0.
\]
Optimisation yields
$\|f\|_{L^{\infty}}\le C\|f\|_{L^{p}}^{1-\sigma}\|\nabla f\|_{L^{m}}^\sigma$.
If the balancing radius exceeds $\ell_0$, use zero-mean
Poincar\'e--Morrey. Combining this with
\[
 \|f\|_{L^{p}}\le C\sqrt p\,\|f\|_{L^{2}}^{2/p}
                              \|\nabla f\|_{L^{2}}^{1-2/p}
\]
proves the claim; see \cite[Lemma~2.2]{HuangLi2016}.
\end{proof}

\begin{proposition}[Uniform weighted and local flux bounds]\label{per:prop:flux}
For every $0<\tau<S$,
\begin{equation}\label{per:eq:weighted-flux}
 \int_\tau^S\frac{\|G(t)\|_{L^{\infty}}^2}{e+A_1^2(t)}\,\mathrm dt
 \le C_\tau R_S^2[\log(e+R_S)]^2,
\end{equation}
where $C_\tau$ is independent of $S$. For $\tau\le s<t\le S$,
\begin{equation}\label{per:eq:local-flux}
 \int_s^t\|G(\xi)\|_{L^{\infty}}\,\mathrm d\xi
 \le C_\tau R_t\log(e+R_t)
           \left(t-s+\int_s^t A_3^2\,\mathrm d\xi\right)^{1/2}.
\end{equation}
\end{proposition}
\begin{proof}
Since $\int\rho\dot u\dd=0$, weighted Poincar\'e and Sobolev imply
$\|\dot u\|_{L^{m}}\le C\sqrt m\,\|\nabla\dot u\|_{L^{2}}$.
Periodic Calder\'on--Zygmund therefore gives
\begin{equation}\label{per:eq:high-gradient}
 \|\nabla G\|_{L^{m}}\le Cm\|\rho\dot u\|_{L^{m}}
             \le Cm^{3/2}R_S\|\nabla\dot u\|_{L^{2}},\qquad m\ge2.
\end{equation}
Apply Lemma~\ref{per:lem:interpolation} to $G-\bar G$, with
\[
 p=8+(b+1)\log(e+R_S),\qquad m=e^p.
\]
Then $\alpha,\sigma\le C/p$, $\sigma-\alpha\le C/(pm)$, and
$R_S^{C\alpha+C\sigma}m^{3\sigma/2}\le C$.
Equations \eqref{eq:flux-elliptic} and \eqref{per:eq:high-gradient} yield
\begin{equation}\label{per:eq:pointwise-flux}
\begin{aligned}
 \frac{\|G-\bar G\|_{L^{\infty}}}{\sqrt{e+A_1^2}}
 \le{}&C\sqrt{pR_S}
 \left(\frac{A_1^2}{e+A_1^2}\right)^{\alpha/2}
 \left(\frac{A_2^2}{e+A_1^2}\right)^{(1-\alpha-\sigma)/2}\\
 &\times\left(\frac{\|\nabla\dot u\|_{L^{2}}}{\sqrt{e+A_1^2}}\right)^\sigma.
\end{aligned}
\end{equation}
Retaining the $A_1^2$ factor supplies global time integrability.
Insert $e+A_2^2$ in the last ratio and split
$(e+A_2^2)^{\sigma/2}\le e^{\sigma/2}+A_2^\sigma$.
For the first term, apply H\"older's inequality with weights
$(\alpha,1-\alpha-\sigma,\sigma)$ to the integrals of
$A_1^2/(e+A_1^2)$, $A_2^2/(e+A_1^2)$, and
$\|\nabla\dot u\|_{L^{2}}^2/(e+A_2^2)$.
For the second term, \eqref{per:eq:log-acceleration} gives
\[
 \sup_{\tau\le t\le S}
 \left(\frac{A_2^2}{e+A_1^2}\right)^{(\sigma-\alpha)/2}
 \le\exp\!\left(
   \frac{C_\tau R_S^{b+1}\log(e+R_S)}{pm}\right)\le C_\tau,
\]
after which H\"older uses weights $(1-\sigma,\sigma)$.
Using \eqref{per:eq:energy-comparison}, \eqref{eq:log-energy}, and
\eqref{per:eq:log-acceleration}, we obtain
\[
 \left(\int_\tau^S
       \frac{\|G-\bar G\|_{L^{\infty}}^2}{e+A_1^2}\,\mathrm dt\right)^{1/2}
 \le C_\tau\sqrt p\,R_S^{1+b\sigma/2}
                        [\log(e+R_S)]^{1/2}
 \le C_\tau R_S\log(e+R_S).
\]
Combining this with
$\int_\tau^S|\bar G|^2/(e+A_1^2)\,\mathrm dt
 \le C\int_\tau^S A_3^2\,\mathrm dt\le C$
proves \eqref{per:eq:weighted-flux}. Cauchy--Schwarz on each prefix $[0,t]$,
using $e+A_1^2\le C(1+A_3^2)$, gives \eqref{per:eq:local-flux}.
\end{proof}

\subsection{The density bound}\label{per:sec:density}

\begin{lemma}[Uniform initial interval]\label{per:lem:initial-layer}
There are $\tau_0\in(0,1]$ and $M_0$, depending only on the data and
coefficient bounds in Theorem~\ref{per:thm:main}, such that
\begin{equation}\label{per:eq:initial-layer}
 \sup_{0\le t\le\min\{S,\tau_0\}}\|\rho(t)\|_{L^{\infty}}\le M_0.
\end{equation}
\end{lemma}
\begin{proof}
Let $M_0=2(1+C_0)$ and work before the first contact with $M_0$,
with $t\le\min\{S,1\}$.
Lemma~\ref{lem:log-energy} gives
$\sup A_1^2+\int A_2^2\le C_*$.
Multiplying \eqref{per:eq:acceleration-differential} by time and applying
Gronwall gives
\[
 \sup_{0<s\le t}sA_2^2(s)
             +\int_0^t s\|\nabla\dot u\|_{L^{2}}^2\,\mathrm ds\le C_*.
\]
The initial weighted boundary term vanishes for each smooth solution.
By \eqref{per:eq:high-gradient} with $m=4$ and Gagliardo--Nirenberg,
\[
 \|G\|_{L^{\infty}}\le C_*\bigl(1+\|\nabla\dot u\|_{L^{2}}^{2/3}\bigr).
\]
Consequently,
\[
\begin{aligned}
 \int_0^t\|G\|_{L^{\infty}}\,\mathrm ds
 &\le C_*t+C_*
 \left(\int_0^t s\|\nabla\dot u\|_{L^{2}}^2\,\mathrm ds\right)^{1/3}
 \left(\int_0^t s^{-1/2}\,\mathrm ds\right)^{2/3}\\
 &\le C_*t^{1/3}.
\end{aligned}
\]
Using the periodic case of \eqref{eq:theta-flow}, we obtain along every trajectory
$\theta(\rho(t))\le\theta(C_0)+C_*t^{1/3}$.
Choose $\tau_0$ so that
$\theta(C_0)+C_*\tau_0^{1/3}<\theta(M_0)$; this excludes a first
contact before $\tau_0$, independently of the density lower bound.
\end{proof}

\begin{lemma}[Density bound from primitive-level costs]
\label{lem:primitive-level-density}
Let $\Omega=\TT$ or $\mathbb R^2$, and suppose that
\eqref{per:eq:primitive-growth} holds. Assume that the flow is complete,
comparison along positive-density trajectories is valid on $[0,S]$,
and spatial maxima are attained. Suppose also that
\[
 \Dt\theta(\rho)+\rho^\gamma\le |G|+C,
 \qquad \int_0^S A_3^2\,\mathrm dt\le C.
\]
Suppose, with constants independent of $S$, that
\[
 \sup_{0\le t\le\min\{S,\tau_0\}}\|\rho(t)\|_{L^\infty}\le M_0
\]
for some $\tau_0>0$, and that
\[
 \int_s^t\|G(\xi)\|_{L^\infty}\,\mathrm d\xi
 \le C_\tau(1+R_t)\log(e+R_t)
       \left(t-s+\int_s^t A_3^2(\xi)\,\mathrm d\xi\right)^{1/2},
 \qquad \tau\le s<t\le S,
\]
for each fixed $\tau>0$. Then $R_S\le C$ independently of $S$.
\end{lemma}
\begin{proof}
There is nothing to prove if $S\le\tau_0$. Otherwise fix
$\tau=\tau_0/2$ and choose
\[
 r_0>\max\{2,r_\lambda,M_0,(2C)^{1/\gamma}\},
 \qquad \theta(r_0)>0,
\]
large enough for the absorption below. Put
$r_j=\theta^{-1}(2^j\theta(r_0))$ and let $t_j$ be the first time
the running maximum reaches $r_j$, with $t_0=\tau$.
For attained levels set
$\mathcal E_j=\int_{t_{j-1}}^{t_j}A_3^2\,\mathrm dt$, so that
$\sum_j\mathcal E_j\le C$.
On a trajectory reaching $r_j$ at $t_j$, let $\sigma_j$ be its last
crossing of $r_{j-1}$. The differential inequality and the flux bound give
\[
\begin{aligned}
 \frac12\theta(r_j)
 &\le C_\tau r_j\log(e+r_j)
          \sqrt{t_j-\sigma_j+\mathcal E_j}
                  -\frac12r_{j-1}^\gamma(t_j-\sigma_j)\\
 &\le C_\tau r_j\log(e+r_j)\sqrt{\mathcal E_j}
       +C_\tau\frac{r_j^2[\log(e+r_j)]^2}{r_{j-1}^\gamma}.
\end{aligned}
\]
Call $j$ good if
$\log r_j\le(\gamma+1)\log r_{j-1}/2$ and bad otherwise.
On good levels the last term is bounded by
$C_\tau r_j^{2/(\gamma+1)}[\log(e+r_j)]^2
\le\theta(r_j)/4$ after increasing $r_0$. Thus
\begin{equation}\label{eq:primitive-level-cost}
 \mathcal E_j\ge
 c\left[\frac{\theta(r_j)}{r_j\log(e+r_j)}\right]^2
 \ge c[\log(e+r_j)]^{2\kappa-2}
 \ge c(j+1)^{2\kappa-2}.
\end{equation}
The last step uses $r_j\le C\theta(r_j)=C2^j\theta(r_0)$ and
$\kappa\le1$. If $j_\ell$ is the $\ell$th bad index, then
\[
 \left(\frac{\gamma+1}{2}\right)^\ell\log r_0
 \le\log r_{j_\ell}\le C(j_\ell+1),
 \qquad
 \sum_{j\ {\rm bad}}\frac1{j+1}<\infty.
\]
Since $2\kappa-2\ge-1$, every attained $N$ satisfies
\[
 C\ge\sum_{j=1}^N\mathcal E_j
 \ge c\sum_{\substack{1\le j\le N\\j\ {\rm good}}}\frac1{j+1}
 \ge c\sum_{j=1}^N\frac1{j+1}-C.
\]
Hence only a bounded number of levels can be attained, and $R_S\le C$.
\end{proof}

\begin{proof}[Proof of Proposition~\ref{per:prop:apriori}]
By \eqref{eq:theta-flow} and \eqref{per:eq:basic},
$\Dt\theta(\rho)+P\le |G|+C$ and $\int_0^S A_3^2\,dt\le C$.
Together with Lemma~\ref{per:lem:initial-layer} and
\eqref{per:eq:local-flux}, these bounds verify the hypotheses of
Lemma~\ref{lem:primitive-level-density}, which yields $R_S\le C$.

Substitution into \eqref{eq:log-energy}, \eqref{per:eq:energy-comparison},
and \eqref{per:eq:basic} gives
\[
 \sup_{0\le t\le S}(A_1^2+A_3^2+\|u\|_{H^1}^2)
                      +\int_0^S A_2^2\,\mathrm dt\le C.
\]
Multiply \eqref{per:eq:acceleration-differential} by $\min\{1,t\}$ and
use $0\le(\min\{1,t\})'\le\mathbf1_{(0,1)}$. Gronwall's inequality gives
\[
 \sup_{0<t\le S}\min\{1,t\}A_2^2(t)
       +\int_0^S\min\{1,t\}\|\nabla\dot u\|_{L^{2}}^2\,\mathrm dt\le C.
\]
The initial weighted boundary term vanishes for each smooth
approximation. Returning to the original frame proves
\eqref{per:eq:main-bound}.
\end{proof}

\section{Whole-space potentials and uniform density estimates}\label{sec:cauchy-density}

Throughout this section, $\Omega=\R^2$ and the common structural
condition~\ref{per:ass:lambda} is imposed. The estimates are first
proved on compact intervals of the local strong solution, with
constants independent of the interval endpoint.

\subsection{Potential estimates without spatial localization}
\label{ws:sec:potentials}

The following estimates replace periodic Poincar\'e estimates without
requiring uniform spatial localization of the mass.
Write $I=(-\Delta)^{-1}$ and $T_{ij}=\partial_i\partial_j I$.

\begin{lemma}[A weighted potential estimate]\label{ws:lem:potential}
Let $0\le\rho\le R$, $\int\rho=M$, $v\in H^1_{\rm loc}$,
$\|\sqrt\rho\,v\|_{L^2}<\infty$, and $\|\nabla v\|_{L^2}<\infty$.
Then
\begin{equation}\label{ws:eq:potential-bound}
 \|\nabla I\cdot(\rho v)\|_{L^\infty}
 \le C\sqrt R\,\|\sqrt\rho\,v\|_{L^2}
       \left[1+\log\left(1+
       \frac{\sqrt M\,\|\nabla v\|_{L^2}}
            {\|\sqrt\rho\,v\|_{L^2}}\right)\right]^{1/2}.
\end{equation}
The right side is interpreted as zero when
$\|\sqrt\rho\,v\|_{L^2}=0$.
\end{lemma}
\begin{proof}
The kernel is bounded by $C|z|^{-1}$. Put $L=\sqrt{M/R}$.
Fix $x$, $0<\ell\le L$, and let $v_\ell$ be the ordinary average
of $v$ over $B_\ell(x)$. If $m=\int_{B_\ell(x)}\rho$, splitting a
radial integral at $\sqrt{m/R}$ gives
\[
 \int_{B_\ell(x)}\frac{\rho(y)}{|x-y|}\,dy\le C\sqrt{Rm},\qquad
 \sqrt m\,|v_\ell|
 \le \|\sqrt\rho\,v\|_{L^2}
       +C\sqrt R\,\ell\|\nabla v\|_{L^2}.
\]
Scaled Poincar\'e and H\"older inequalities give
\[
 \int_{B_\ell(x)}
       \frac{\rho|v-v_\ell|}{|x-y|}\,dy
 \le CR\ell\|\nabla v\|_{L^2}.
\]
The near part is therefore bounded by
$C\sqrt R\,\|\sqrt\rho\,v\|_{L^2}+CR\ell\|\nabla v\|_{L^2}$.
For the far part, use weighted Cauchy--Schwarz and
\[
 \int_{|x-y|\ge\ell}\frac{\rho(y)}{|x-y|^2}\,dy
 \le CR\left(1+\log\frac L\ell\right)
\]
to obtain $C\sqrt R\,\|\sqrt\rho\,v\|_{L^2}[1+\log(L/\ell)]^{1/2}$.
Choosing
\[
 \ell=\frac{L\|\sqrt\rho\,v\|_{L^2}}
                 {\|\sqrt\rho\,v\|_{L^2}
                       +\sqrt M\,\|\nabla v\|_{L^2}}
\]
proves the claim when $\|\sqrt\rho v\|_{L^2}>0$; otherwise $\rho v=0$.
These bounds on the absolute kernel justify the representation by truncation.
\end{proof}

\begin{lemma}[A potential interpolation estimate]\label{ws:lem:potential-interpolation}
Let $v\in H^1_{\rm loc}$,
$h=\nabla I\cdot(\rho v)\in L^2$, and $0\le\rho\le R$.
Assume $\|\sqrt\rho v\|_{L^2}<\infty$ and
$\|\nabla v\|_{L^2}<\infty$.
For every $p\ge4$, with $\sigma=2/(p+2)$,
\begin{equation}\label{ws:eq:potential-interpolation}
 \|h\|_{L^\infty}
 \le C\sqrt R\,\|\sqrt\rho v\|_{L^2}+
 C\sqrt{pR}\,\|h\|_{L^2}^\sigma
      \|\sqrt\rho v\|_{L^2}^{1-2\sigma}
      \|\nabla v\|_{L^2}^\sigma .
\end{equation}
The constant is independent of $p$.
\end{lemma}
\begin{proof}
Let $\eta_\ell$ be a nonnegative smooth averaging kernel of radius
$\ell$, and let $K$ be the kernel of $\nabla I$.
On $|z|<2\ell$, apply the preceding near-field estimate to $K$. The bound
$|\eta_\ell*K|\le C/\ell$ contributes at most
$C\sqrt R\,\|\sqrt\rho v\|_{L^2}$.
On the complement,
\[
 |K-\eta_\ell*K|(z)\le C\ell|z|^{-2},\qquad
 \int_{|z|>2\ell}\rho(x-z)|z|^{-4}\,dz\le CR/\ell^2.
\]
Thus weighted Cauchy--Schwarz gives
\[
 \|h-\eta_\ell*h\|_{L^\infty}
 \le C\sqrt R\,\|\sqrt\rho v\|_{L^2}
       +CR\ell\|\nabla v\|_{L^2}.
\]
Also $\|\nabla h\|_{L^2}\le C\sqrt R\,\|\sqrt\rho v\|_{L^2}$.
The quantitative Gagliardo--Nirenberg inequality
\cite[Lemma~2.2]{HuangLi2022} gives
\[
 \|\eta_\ell*h\|_{L^\infty}
 \le C\ell^{-2/p}\|h\|_{L^p}
 \le C\sqrt p\,\ell^{-2/p}\|h\|_{L^2}^{2/p}
                 (\sqrt R\,\|\sqrt\rho v\|_{L^2})^{1-2/p}.
\]
Balancing the last term with $R\ell\|\nabla v\|_{L^2}$ gives
\[
 C\sqrt R\,(\sqrt p)^{p/(p+2)}
       \|h\|_{L^2}^{2/(p+2)}
       \|\sqrt\rho v\|_{L^2}^{(p-2)/(p+2)}
       \|\nabla v\|_{L^2}^{2/(p+2)}.
\]
Since $(\sqrt p)^{p/(p+2)}\le\sqrt p$, this proves
\eqref{ws:eq:potential-interpolation}; zero factors follow by a limit.
\end{proof}

\subsection{Basic energy and spatial moments}

Work on a compact interval $[0,T]$ of a local strong solution.
The constants in the estimates below are independent of $T$, the provisional
density upper bound, and a positive density lower bound. To prove mass
conservation, test continuity with expanding cutoffs. The resulting boundary term is bounded by
$C L^{-1}\|\sqrt\rho u\|_{L^2}\|\rho\|_{L^1}^{1/2}$. After integration over a
finite time interval, this term vanishes as $L\to\infty$; the local strong
norms justify the limit. Similarly, momentum conservation follows by testing the conservative
momentum equation with logarithmic spatial cutoffs whose gradients tend
to zero in $L^2$ and $L^\infty$. The local strong class gives stress and pressure in $L^2$ and
$\rho|u|^2$ in $L^1$ on every compact existence interval. Hence all boundary terms vanish.

For the low-order estimates, we use the Galilean frame determined by
the conserved momentum, so that $\int\rho u\dd=0$. Density bounds, velocity gradients,
the effective flux, and material acceleration are unchanged by this
transformation.
The energy balance is
\begin{equation}\label{ws:eq:physical-energy-balance}
 \frac{d}{dt}\int\left(\frac12\rho|u|^2+
                      \frac{\rho^\gamma}{\gamma-1}\right)\dd
       +A_3^2=0.
\end{equation}
For the cutoff argument, temporarily set
$R=1+\sup_{t\le T}\|\rho(t)\|_{L^\infty}$, which is finite by local strong regularity. Choose $\chi_L=1$
on $B_L$, vanishing outside $B_{L^2}$, with
$\|\nabla\chi_L\|_{L^2}\le C/\sqrt{\log L}$ and
$\|\nabla\chi_L\|_{L^\infty}\to0$.
Test momentum by $\chi_L DF_N(u)$, where
\[
 \begin{gathered}
 F_N(z)=N^2\left(\sqrt{1+|z|^2/N^2}-1\right),
 \qquad |DF_N|\le N,\\
 0\le D^2F_N\le I_2,\qquad
 0\le F_N(z)\le\min\{|z|^2/2,N|z|\}.
 \end{gathered}
\]
Use continuity for the truncated kinetic energy $\chi_L\rho F_N(u)$.
For fixed $N$, the convective and stress--pressure boundary terms are
bounded, respectively, by
\[
 CN\|\nabla\chi_L\|_{L^\infty}\int\rho|u|^2,\qquad
 CN\|\nabla\chi_L\|_{L^2}
       (\|P\|_{L^2}+C_R\|\nabla u\|_{L^2}),
\]
and vanish as $L\to\infty$ after time integration.
Dominated convergence then permits $N\to\infty$, with majorants
$C_R|\nabla u|^2+CP|\nabla u|$.
For the internal-energy equation, use
$\int P|u|\le R^{\gamma-1}\sqrt M\|\sqrt\rho u\|_{L^2}$.
Adding these identities and applying Hodge yields
\eqref{ws:eq:physical-energy-balance}. The temporary $R$ is used only to
justify the calculation. The identity extends to time zero by the initial
kinetic trace in Lemma~\ref{ws:lem:general-local}, justified by
\eqref{ws:eq:local-kinetic-tail}. Hence
\begin{equation}\label{ws:eq:basic-energy}
 \sup_{t\le T}\int(\rho|u|^2+\rho^\gamma)\dd+
        \int_0^TA_3^2(t)\,dt\le C.
\end{equation}
In particular $\int_0^T\|\nabla u\|_{L^2}^2\,dt\le C$.
Since $a<2$, $|\nabla\bar x^{a/2}|\le C$. Continuity and
\eqref{ws:eq:basic-energy} give
\[
 \frac{d}{dt}\left(\int\rho\bar x^a\dd\right)^{1/2}
       \le C\|\sqrt\rho u\|_{L^2}\le C.
\]
Consequently,
\begin{equation}\label{ws:eq:finite-weighted-mass}
 \int\rho\bar x^a\dd\le C(e+t)^2,\qquad
 \int\rho\log(e+|x|)\dd\le C\log(e+t),
\end{equation}
where the second assertion follows from Jensen's inequality.
On every interval $[0,S]$ with $S\le T$, choosing $N_T$ sufficiently large gives
\begin{equation}\label{ws:eq:finite-mass-localization}
 \inf_{t\le S}\int_{B_{N_T}}\rho\dd\ge M/4.
\end{equation}
These spatial estimates use only the basic energy and mass conservation.

\begin{lemma}[A weighted Poincar\'e inequality]\label{ws:lem:log-weighted-poincare}
Let $0\le\rho\in L^1\cap L^\gamma$, $\gamma>1$, have mass $M>0$ and finite
logarithmic spatial moment. If $\nabla v\in L^2$ and $\sqrt\rho v\in L^2$, then
\begin{equation}\label{ws:eq:log-weighted-poincare}
 \int\rho|v-v_\rho|^2\dd
 \le C_\gamma\left(\|\rho\|_{L^\gamma}
                  +\int\rho\log(e+|x|)\dd\right)
                   \|\nabla v\|_{L^2}^2,
 \qquad v_\rho=M^{-1}\int\rho v\dd.
\end{equation}
\end{lemma}
\begin{proof}
Let $v_n$ be the mean on $Q_n=n+[0,1]^2$.
H\"older and local Sobolev--Poincar\'e give
\[
 \sum_n\int_{Q_n}\rho|v-v_n|^2\dd
 \le C_\gamma\|\rho\|_{L^\gamma}\|\nabla v\|_{L^2}^2.
\]
Homogeneous Sobolev duality gives
\[
 |v_n-v_0|^2
 \le\|\mathbf1_{Q_n}-\mathbf1_{Q_0}\|_{\dot H^{-1}}^2
                   \|\nabla v\|_{L^2}^2
 \le C\log(e+|n|)\|\nabla v\|_{L^2}^2.
\]
To obtain the last bound, split the Fourier integral at
$|\xi|=|n|^{-1}$ and $1$; bounded values of $n$ are treated directly.
Summing these bounds against the cell masses and minimizing over constants
gives \eqref{ws:eq:log-weighted-poincare}. Approximation modulo constants
extends the estimate to the homogeneous energy class.
\end{proof}

In the zero-momentum frame, \eqref{ws:eq:log-weighted-poincare} and
\eqref{ws:eq:finite-weighted-mass} imply
\begin{equation}\label{ws:eq:mean-control-before-density}
 \|\sqrt\rho u\|_{L^2}^2+|\bar u_{Q_0}|^2
 \le C\log(e+t)\|\nabla u\|_{L^2}^2.
\end{equation}
Here $\bar u_{Q_0}$ is the mean on $Q_0$. For its bound, use
$M|\bar u_{Q_0}|^2\le\int\rho|u-\bar u_{Q_0}|^2\dd$
and the cell estimates in the preceding proof.

\subsection{Density moments and preliminary energy decay}

\begin{lemma}[Subpolynomial density moments]\label{ws:lem:finite-moment}
Assume the structural condition~\ref{per:ass:lambda} and
\eqref{ws:eq:initial-data}. On the interval of existence of the local strong
solution, every integer $k\ge3$ with $(2\kappa-1)k\ge1$ satisfies
\begin{equation}\label{ws:eq:subpolynomial-moments}
 \int\left(\rho\theta(\rho)_+^k+
       \rho^{k+1}[\log(e+\rho)]^{\kappa k}\right)\dd
 \le C_k\exp\!\bigl(C_k\sqrt{\log(e+t)}\bigr).
\end{equation}
The constants depend only on $k$, the low-order initial data, and
the coefficient assumptions, independently of the provisional density bound.
\end{lemma}
\begin{proof}
For a sufficiently small fixed $c>0$, we first prove
\begin{equation}\label{ws:eq:weighted-trudinger}
 \int(e+t)^{-1}\bar x^{-4}
       \exp\!\left(\frac{c|u|^2}{\|\nabla u\|_{L^2}^2}\right)\dd
 \le C.
\end{equation}
Indeed, let $\bar u_{Q_n}$ be the unit-cube means. The proof of
Lemma~\ref{ws:lem:log-weighted-poincare} and
\eqref{ws:eq:mean-control-before-density} give
\[
 |\bar u_{Q_n}|^2
 \le C\bigl(\log(e+t)+\log(e+|n|)\bigr)\|\nabla u\|_{L^2}^2.
\]
The mean-zero Trudinger inequality on each unit cube gives, for sufficiently small fixed $c$,
\[
 \int_{Q_n}\exp\!\left(\frac{c|u|^2}{\|\nabla u\|_{L^2}^2}\right)\dd
 \le C(e+t)^{C c}(e+|n|)^{C c}.
\]
Choose $c$ so that $Cc<1$ and sum against $\bar x^{-4}$.
If $\nabla u=0$, then $u=0$ in the zero-momentum frame, and we interpret
the assertion accordingly.

Put
\[
 \psi=(-\Delta)^{-1}\diver(\rho u),\qquad
 \varphi=(\theta(\rho)-\psi)_+,\qquad
 \mathcal M_k=\int\rho\varphi^k\dd.
\]
By \eqref{eq:potential-comparison} and
$\theta(\rho)_+\le\varphi+|\psi|$, absorption gives
\begin{equation}\label{ws:eq:finite-moment-comparison}
 \int\left(\rho\theta(\rho)_+^k+
       \rho^{k+1}[\log(e+\rho)]^{\kappa k}\right)\dd
 \le C_k(1+\mathcal M_k).
\end{equation}

Orlicz--Young with respect to $(e+t)^{-1}\bar x^{-4}\dd$ and
\eqref{ws:eq:weighted-trudinger} give, with $p=k+1$,
\[
 \int\rho^p|u|^p\dd
 \le C_k\|\nabla u\|_{L^2}^p
 \left[1+\int\rho^p
 \bigl(\log(e+\rho)+\log(e+t)+\log\bar x\bigr)^{p/2}\dd\right].
\]
To control the spatial logarithm, put $Y=1+\bar x/(e+t)^{2/a}$. Then
\[
 \log\bar x\le\frac2a\log(e+t)+\log Y,
 \qquad \int\rho Y^a\dd\le C
\]
by \eqref{ws:eq:finite-weighted-mass}. Splitting at
$\log Y=(4k/a)\log(e+\rho)$ gives
\[
 \rho^p(\log Y)^{p/2}
 \le C_k\rho^p[\log(e+\rho)]^{p/2}+C_k\rho Y^a.
\]
Since $p/2\le\kappa k$, \eqref{ws:eq:finite-moment-comparison} implies
\begin{equation}\label{ws:eq:log-momentum-moment}
 \|\rho u\|_{L^{k+1}}
 \le C_k\sqrt{\log(e+t)}\,\|\nabla u\|_{L^2}
                         (1+\mathcal M_k)^{1/(k+1)}.
\end{equation}
The time-dependent measure is used only in this spatial estimate and
is never differentiated.

Momentum and continuity give
\[
 \Dt(\theta(\rho)-\psi)+P
       =-\sum_{i,j}[u_i,T_{ij}](\rho u_j).
\]
Equations~\eqref{eq:moment-commutator} and
\eqref{ws:eq:log-momentum-moment} give
\begin{equation}\label{ws:eq:log-time-moment-ode}
 \mathcal M_k'+k\int\rho^{\gamma+1}\varphi^{k-1}\dd
 \le C_k\sqrt{\log(e+t)}\,A_3^2(1+\mathcal M_k).
\end{equation}
On each compact existence interval, $\psi$ is bounded and $\rho$ tends
uniformly to zero at spatial infinity. Thus $\varphi$ has bounded support
and vanishes for all sufficiently small densities. These cutoffs allow the
transport chain rule to be applied across vacuum without affecting the constants.
The initial density and kinetic energy give $\mathcal M_k(0)\le C_k$, and
the basic energy implies
\[
 \int_0^t\sqrt{\log(e+s)}\,A_3^2(s)\,ds
 \le C\sqrt{\log(e+t)}.
\]
Gronwall and \eqref{ws:eq:finite-moment-comparison} prove
\eqref{ws:eq:subpolynomial-moments}.
\end{proof}

\begin{lemma}[Preliminary decay and weighted dissipation]
\label{ws:lem:preliminary-energy}
Under the structural condition~\ref{per:ass:lambda} and
\eqref{ws:eq:subpolynomial-moments}, set
\[
 \mathcal E(t)=\int\left(
       \frac12\rho|u-u_\infty|^2+\frac{P}{\gamma-1}\right)\dd.
\]
On every existence interval $[0,S]$,
\begin{equation}\label{ws:eq:preliminary-energy-decay}
 \mathcal E(t)\le C
       \left[\frac{\log(e+t)}{\sqrt{1+t}}\right]^{\frac{\gamma-1}\gamma},
       \qquad 0\le t\le S,
\end{equation}
and
\begin{equation}\label{ws:eq:log-weighted-dissipation}
 \int_0^S\sqrt{\log(e+t)}\,A_3^2(t)\,dt\le C.
\end{equation}
The constants are independent of $S$ and the density upper bound.
\end{lemma}
\begin{proof}
We continue in the zero-momentum frame, where
$\mathcal E'=-A_3^2$ and $\mathcal E\le\mathcal E(0)$.

Write $\psi=I\diver(\rho u)$. For a nonnegative
$g\in C^1([0,\infty))$ satisfying $g(r)=O(r^\gamma)$ at zero, set
\[
 h(r)=rg'(r)-g(r),\qquad
 \mathcal H(r)=r\int_0^r
       \frac{(2\mu+\lambda(z))g(z)}{z^2}\,dz,\qquad
 J=\int\rho u\cdot\nabla I g(\rho)\dd.
\]
Since $g(\rho)_t+\diver(g(\rho)u)=-h(\rho)\diver u$ and
$r\mathcal H'-\mathcal H=(2\mu+\lambda)g$, testing momentum by
$\nabla I g(\rho)$ and using continuity gives
\begin{equation}\label{ws:eq:preliminary-pressure-test}
 \left(J+\int\mathcal H(\rho)\dd\right)'
       +\int P g(\rho)\dd
 =\sum_{i,j}\int\rho u_i[u_j,T_{ij}]g(\rho)\dd
       +\int\psi h(\rho)\diver u\dd.
\end{equation}
Spatial and kernel regularization justify this identity on compact
existence intervals. The condition at zero ensures that $\mathcal H(\rho)$
is integrable; the estimates below apply to both choices of $g$ used here.

Choose
\[
 g(r)=\frac{r^\gamma}{2\mu+\lambda(r)+r^\gamma}.
\]
Then
\[
\begin{gathered}
 0\le g(r)\le\min\{1,r^\gamma/(2\mu)\},\qquad
 |h(r)|\le Cg(r),\\
 (2\mu+\lambda(r))g(r)^2\le r^\gamma/4,
 \qquad 0\le\mathcal H(r)\le\frac{r^\gamma}{\gamma-1}.
\end{gathered}
\]
The elementary kernel bound gives
\[
\begin{aligned}
 \|\nabla I g(\rho)\|_{L^\infty}
 &\le C\|g(\rho)\|_{L^1}^{1/2}\|g(\rho)\|_{L^\infty}^{1/2}
       \le C\|P\|_{L^1}^{1/2},\\
 |J|&\le C\|\sqrt\rho u\|_{L^2}\|P\|_{L^1}^{1/2}
       \le C\mathcal E.
\end{aligned}
\]
Hardy--Littlewood--Sobolev and H\"older give
\[
\begin{aligned}
 \|\rho u\|_{L^{2\gamma/(\gamma+1)}}+\|\psi\|_{L^{2\gamma}}
 &\le C\|\sqrt\rho u\|_{L^2}\|P\|_{L^1}^{1/(2\gamma)},\\
 \|g(\rho)\|_{L^{2\gamma/(\gamma-1)}}
       +\|h(\rho)\|_{L^{2\gamma/(\gamma-1)}}
 &\le C\left(\int Pg(\rho)\dd\right)^{(\gamma-1)/(2\gamma)}.
\end{aligned}
\]
The second estimate uses $g^{2\gamma/(\gamma-1)}\le CPg$.
The commutator theorem~\cite{CRW1976},
$\|u\|_{\mathrm{BMO}}\le C\|\nabla u\|_{L^2}$, and
\eqref{ws:eq:mean-control-before-density} bound the right-hand side of
\eqref{ws:eq:preliminary-pressure-test} by
\[
 C\|\sqrt\rho u\|_{L^2}A_3
       \left(\int Pg\dd\right)^{(\gamma-1)/(2\gamma)}
 \le C\sqrt{\log(e+t)}A_3^2,
\]
since $\int Pg\le\int P\le C$.
Choose $\varepsilon>0$ small and set
\[
 \delta(t)=\frac{\varepsilon}{\sqrt{\log(e+t)}},\qquad
 \mathcal F=\mathcal E+\delta\left(J+\int\mathcal H(\rho)\dd\right).
\]
Then $\mathcal F\asymp\mathcal E$ and
\begin{equation}\label{ws:eq:preliminary-energy-differential}
 \mathcal F'+cA_3^2+\delta\int Pg\dd
                  \le C|\delta'|\mathcal F,
 \qquad \int_0^\infty|\delta'|\,dt<\infty.
\end{equation}

H\"older gives
\[
 \|P\|_{L^1}
 \le\left(\int Pg(\rho)\dd\right)^{\frac{\gamma-1}{2\gamma-1}}
 \left(\int\rho(2\mu+\lambda(\rho)+P)^{\frac{\gamma-1}{\gamma}}\dd
       \right)^{\frac\gamma{2\gamma-1}}.
\]
By \eqref{eq:coefficient-bounds} and \eqref{ws:eq:subpolynomial-moments},
\[
\begin{aligned}
 \int\rho(2\mu+\lambda(\rho)+P)^{\frac{\gamma-1}{\gamma}}\dd
 &\le C\int\left(\rho+
      \rho^{1+\frac{b(\gamma-1)}\gamma}+\rho^\gamma\right)\dd\\
 &\le C\exp\!\bigl(C\sqrt{\log(e+t)}\bigr).
\end{aligned}
\]
Consequently,
\[
 \int Pg(\rho)\dd
 \ge c\exp\!\bigl(-C\sqrt{\log(e+t)}\bigr)
              \|P\|_{L^1}^{\frac{2\gamma-1}{\gamma-1}}
 \ge\frac{c}{\sqrt{1+t}}
              \|P\|_{L^1}^{\frac{2\gamma-1}{\gamma-1}}.
\]
Since $\mathcal F\asymp\mathcal E\le C$, the kinetic estimate
\eqref{ws:eq:mean-control-before-density} implies
\[
 \mathcal F^{\frac{2\gamma-1}{\gamma-1}}
 \le C\left(\log(e+t)A_3^2
                  +\sqrt{1+t}\int Pg(\rho)\dd\right).
\]
Thus \eqref{ws:eq:preliminary-energy-differential} yields
\[
 \mathcal F'
 +\frac{c\mathcal F^{\frac{2\gamma-1}{\gamma-1}}}
        {\sqrt{1+t}\log(e+t)}
 \le C|\delta'|\mathcal F.
\]
A bounded integrating factor removes the right-hand side. Since
\[
 \int_0^t\frac{ds}{\sqrt{1+s}\log(e+s)}
       \ge\frac{c\sqrt{1+t}}{\log(e+t)},\qquad t\ge1,
\]
Bernoulli integration proves \eqref{ws:eq:preliminary-energy-decay}.

Finally, integrate $\mathcal E'=-A_3^2$ with weight
$\sqrt{\log(e+t)}$. Dropping the nonpositive terminal term gives
\[
 \int_0^S\sqrt{\log(e+t)}\,A_3^2\,dt
 \le\mathcal E(0)+\frac12\int_0^S
       \frac{\mathcal E(t)}{(e+t)\sqrt{\log(e+t)}}\,dt\le C
\]
by \eqref{ws:eq:preliminary-energy-decay}.
\end{proof}

\begin{lemma}[Uniform density moments]\label{ws:lem:uniform-moments}
Under the assumptions of Lemma~\ref{ws:lem:finite-moment}, let $[0,S]$
be any interval of existence. Every integer $k\ge3$ with
$(2\kappa-1)k\ge1$ satisfies
\begin{equation}\label{ws:eq:uniform-density-moments}
 \sup_{0\le t\le S}\int\left(\rho\theta(\rho)_+^k+
        \rho^{k+1}[\log(e+\rho)]^{\kappa k}\right)\dd\le C_k.
\end{equation}
The constant is independent of $S$ and the provisional density upper bound.
Thus every fixed finite density moment is bounded uniformly throughout the lifespan.
\end{lemma}
\begin{proof}
Lemma~\ref{ws:lem:preliminary-energy} gives
$\int_0^S\sqrt{\log(e+t)}\,A_3^2(t)\,dt\le C$, independently of $S$.
Gronwall in \eqref{ws:eq:log-time-moment-ode} and
\eqref{ws:eq:finite-moment-comparison} prove the claim.
\end{proof}

\subsection{Global pressure dissipation}

Apply \eqref{ws:eq:preliminary-pressure-test} with $g(r)=r^\gamma$ and the
corresponding $J,\mathcal H$, so that $h(\rho)=(\gamma-1)P$.
Lemma~\ref{ws:lem:uniform-moments} bounds every
fixed pressure norm. The same commutator and potential estimates give
\[
\begin{aligned}
 &\left|\sum_{i,j}\int\rho u_i[u_j,T_{ij}]P\dd\right|
       +(\gamma-1)\left|\int\psi P\diver u\dd\right|\\
 &\quad\le C\bigl(\|\rho u\|_{L^{2\gamma/(\gamma+1)}}
                       +\|\psi\|_{L^{2\gamma}}\bigr)
              \|P\|_{L^{2\gamma/(\gamma-1)}}\|\nabla u\|_{L^2}\\
 &\quad\le C\|\sqrt\rho u\|_{L^2}A_3
       \le C\sqrt{\log(e+t)}A_3^2.
\end{aligned}
\]
Moreover,
\[
 \|\nabla IP\|_{L^\infty}
       \le C(\|P\|_{L^1}+\|P\|_{L^4})\le C,
 \qquad |J(t)|\le C\sqrt M\|\sqrt\rho u\|_{L^2}\le C.
\]
For this choice of $g$, $\mathcal H\ge0$ and its initial integral is
finite by \eqref{eq:coefficient-bounds} and \eqref{ws:eq:initial-data}.
Integrating \eqref{ws:eq:preliminary-pressure-test}, dropping the
nonnegative terminal entropy, and using
\eqref{ws:eq:log-weighted-dissipation} give
\begin{equation}\label{ws:eq:pressure-integral}
 \int_0^S\|P(t)\|_{L^2}^2\,dt
 \le C+C\int_0^S\sqrt{\log(e+t)}A_3^2(t)\,dt\le C.
\end{equation}

\subsection{Uniform logarithmic energy and flux estimates}

\begin{lemma}[Uniform weighted and local flux bounds]\label{ws:lem:uniform-flux}
On every existence interval $[0,S]$, put $R=1+R_S$. For $0<\tau<S$,
\begin{equation}\label{ws:eq:weighted-flux}
 \int_\tau^S\frac{\|G\|_{L^\infty}^2}{e+A_1^2}\,dt
 \le C_\tau R^2[\log(e+R)]^2,
\end{equation}
and
\begin{equation}\label{ws:eq:local-flux}
 \int_s^t\|G(\xi)\|_{L^\infty}\,d\xi
 \le C_\tau(1+R_t)\log(e+R_t)
       \left(t-s+\int_s^t A_3^2(\xi)\,d\xi\right)^{1/2},
 \qquad \tau\le s<t\le S.
\end{equation}
The constants are independent of $S$ and the density lower bound.
\end{lemma}
\begin{proof}
Lemma~\ref{ws:lem:uniform-moments} and \eqref{ws:eq:pressure-integral}
give
\[
\begin{aligned}
 A_1^2&\le C(A_3^2+\|P\|_{L^2}^2),&
 A_3^2&\le C(e+A_1^2),\\
 e+A_1^2&\le C(1+A_3^2),&
 \int_0^S A_1^2\,dt&\le C.
\end{aligned}
\]
Continuity of the density and \eqref{ws:eq:initial-gradient-energy} give the
initial upper trace of $A_1$. Thus Lemma~\ref{lem:log-energy} applies
with $R=1+R_S$.
The representation $G=-I\diver(\rho\dot u)$ has no additive harmonic term:
$G\in L^2$, $\rho\dot u\in L^1$, and both expressions have vanishing averages
on large balls.

Bounding the stress square in \eqref{eq:material-acceleration} by
$C(1+\lambda)^2|\nabla u|^2+CP^2$, then using
\eqref{eq:coefficient-bounds} and $P\le R^\gamma$, gives
\begin{equation}\label{ws:eq:hoff-core}
 (A_2^2)'+c\|\nabla\dot u\|_{L^2}^2
 \le C(1+R^{2b})\|\nabla u\|_{L^4}^4
                              +CR^{2\gamma}\|\nabla u\|_{L^2}^2.
\end{equation}
Splitting $|\diver u|\le2P/(2\mu+\lambda(\rho))$ and its complement gives
\[
 \|\diver u\|_{L^4}^4
 \le C\|G\|_{L^4}^4
                   +CR^{2\gamma}\|\nabla u\|_{L^2}^2.
\]
Using Hodge and Gagliardo--Nirenberg,
\begin{equation}\label{ws:eq:fourth-core}
 \|\nabla u\|_{L^4}^4
 \le CR^{b+1}A_1^2A_2^2
                   +CR^{2\gamma}\|\nabla u\|_{L^2}^2.
\end{equation}

There are $\tau_0>0$ and $M_0$ depending only on the low-order
data such that
\begin{equation}\label{ws:eq:initial-layer}
 \sup_{t\le\min(S,\tau_0)}\|\rho(t)\|_{L^\infty}\le M_0 .
\end{equation}
Before first contact with a fixed
$M_0>2(1+\|\rho_0\|_{L^\infty})$ and for $t\le1$,
\eqref{eq:log-energy} bounds $\sup A_1^2+\int A_2^2$ by a
data constant. Multiply \eqref{ws:eq:hoff-core} by time and use
\eqref{ws:eq:fourth-core} to obtain
\[
 \sup tA_2^2+\int t\|\nabla\dot u\|_{L^2}^2\,dt\le C .
\]
Lemma~\ref{ws:lem:potential-interpolation} with $p=4$ gives
$\|G\|_{L^\infty}\le C(A_2+A_2^{1/3}\|\nabla\dot u\|_{L^2}^{1/3})$.
H\"older's inequality, with weights $1/6,1/6,2/3$ on
$A_2^2,t\|\nabla\dot u\|_{L^2}^2,t^{-1/4}$, yields
$\int_0^t\|G\|_{L^\infty}\le C\sqrt t$.
Equation~\eqref{eq:theta-flow} excludes first contact before a fixed
$\tau_0>0$. The initial weighted boundary term vanishes on each smooth
approximation, so the constants are uniform.

It suffices to take $0<\tau<\tau_0$, since larger values follow by restriction. Divide
\eqref{ws:eq:hoff-core} by $e+A_2^2$ and use
\eqref{ws:eq:fourth-core}, \eqref{ws:eq:pressure-integral}, and the
initial-layer bound at $\tau$.
The right-hand side after division is bounded by
\[
 CR^{3b+1}A_1^2+CR^{2b+2\gamma}\|\nabla u\|_{L^2}^2.
\]
Since \eqref{ws:eq:pressure-integral} gives
$\int_0^S A_1^2\,dt\le C$, the fixed choice
$k=3b+2\gamma+4$ suffices in
\begin{equation}\label{ws:eq:acceleration-polynomial}
\begin{aligned}
 &\int_\tau^S
   \frac{\|\nabla\dot u\|_{L^2}^2}
        {e+A_2^2}\,dt\\
 &\qquad+\sup_{\tau\le t\le S}
      \log\bigl(e+A_2^2(t)\bigr)
       \le C_\tau R^k .
\end{aligned}
\end{equation}
Apply Lemma~\ref{ws:lem:potential-interpolation} to
$G=-I\diver(\rho\dot u)$. With $\sigma=2/(p+2)$, it gives
\[
\begin{aligned}
 \frac{\|G\|_{L^\infty}^2}{e+A_1^2}
 &\le CR\frac{A_2^2}{e+A_1^2}
   +CpR^{1+b\sigma}
       \left(\frac{A_1^2}{e+A_1^2}\right)^\sigma
       \left(\frac{A_2^2}{e+A_1^2}\right)^{1-2\sigma}\\
 &\qquad\times
       \left(\frac{\|\nabla\dot u\|_{L^2}^2}{e+A_2^2}\right)^\sigma
       \left(\frac{e+A_2^2}{e+A_1^2}\right)^\sigma.
\end{aligned}
\]
Split $(e+A_2^2)^\sigma$ into its constant and acceleration parts.
H\"older's inequality, \eqref{eq:log-energy}, the acceleration
bound, and $\int_\tau^S A_1^2/(e+A_1^2)\,dt\le C$ yield
\[
\begin{aligned}
 \int_\tau^S\frac{\|G\|_{L^\infty}^2}{e+A_1^2}\,dt
 &\le C_\tau R^2\log(e+R)\\
 &\quad+C_\tau pR^{1+(b+2k)\sigma}
    \bigl([R\log(e+R)]^{1-2\sigma}
          +[R\log(e+R)]^{1-\sigma}\bigr).
\end{aligned}
\]
Taking $p=8+4(b+2k+2)\log(e+R)$ yields
\eqref{ws:eq:weighted-flux}.
Applying the bound on each prefix, with
$e+A_1^2\le C(1+A_3^2)$, Cauchy--Schwarz gives
\eqref{ws:eq:local-flux}.
\end{proof}

\subsection{The time-independent density bound}\label{ws:sec:density}

\begin{proposition}[Uniform whole-space estimates]\label{ws:prop:apriori}
Under the structural condition~\ref{per:ass:lambda} and
\eqref{ws:eq:initial-data}, the local strong solution satisfies
\eqref{ws:eq:uniform-core} and \eqref{ws:eq:flux-square-global}, with all
suprema and integrals restricted to its lifespan. The constants are
independent of the lifespan endpoint and of the density lower bound.
\end{proposition}
\begin{proof}
The local strong class gives at most linear spatial growth of $u$, with
locally time-integrable coefficients, because
$\bar x^{-1}u\in L^\infty_tL^2_x$ and
$\nabla u\in L^1_tL^\infty_x$. Thus the forward and backward flows are
complete. Also $\rho(t)\in L^1\cap W^{1,q}$, $q>2$, tends to zero at
infinity, so its maximum is attained. Vacuum paths remain at zero.
Equations~\eqref{eq:theta-flow}, \eqref{ws:eq:basic-energy},
\eqref{ws:eq:initial-layer}, and \eqref{ws:eq:local-flux} now verify
the assumptions of Lemma~\ref{lem:primitive-level-density}. Hence
$\sup_{t\le S}\|\rho(t)\|_{L^\infty}\le C$.

Substitution in \eqref{eq:log-energy} gives
$\sup A_1^2+\int A_2^2\le C$; Hodge and the fixed pressure moments
give $\sup\|\nabla u\|_{L^2}^2\le C$.
Multiplying \eqref{ws:eq:hoff-core} by $\min\{1,t\}$ and using
\eqref{ws:eq:fourth-core} gives the remaining time-weighted estimates
by Gronwall, since its coefficient is $CA_1^2$ and its source is integrable.
Use $\min\{(t-\varepsilon)_+,1\}$ and let $\varepsilon\downarrow0$
to avoid an initial acceleration trace.
The basic energy and \eqref{ws:eq:pressure-integral} give the other
integrals in \eqref{ws:eq:uniform-core}. Finally,
\eqref{ws:eq:weighted-flux} and $\sup A_1^2\le C$ give
\eqref{ws:eq:flux-square-global}. All estimates are independent of $S$.
Returning to the original Galilean frame proves the proposition.
\end{proof}

\section{Construction and uniqueness of global strong solutions}\label{sec:strong-solutions}

\subsection{Common strong estimates}

\begin{lemma}[Integrability at the initial time]\label{lem:force-integrability}
Let $\Omega=\TT$ or $\R^2$, $q>2$, and $0<S\le T$.
Let $(\rho,u)$ be a local strong solution satisfying the low-order estimates
\eqref{per:eq:main-bound} on $\TT$, or
\eqref{ws:eq:uniform-core}, \eqref{ws:eq:finite-weighted-mass}, and
\eqref{ws:eq:finite-mass-localization} on $\R^2$, restricted to $[0,S]$.
For every $1<s<q/(q-1)$,
\begin{equation}\label{eq:force-integrability}
 \int_0^S\left(\|\rho\dot u\|_{L^q}^s
                   +t\|\rho\dot u\|_{L^q}^2\right)\,dt
       \le C_{q,s,T}.
\end{equation}
On $\TT$, also $\int_0^S t\|\dot u\|_{H^1}^2\,\mathrm dt\le C_T$.
The constants are independent of $S$ and of a density lower bound.
\end{lemma}
\begin{proof}
Weighted Poincar\'e on $\TT$ and the weighted Hardy inequality
\cite[Lemma~2.4]{LiLiang2014} on $\R^2$ give, for every finite $r>2$,
\[
 \|\rho v\|_{L^r}
 \le C_{T,r}\bigl(\|\sqrt\rho v\|_{L^2}+\|\nabla v\|_{L^2}\bigr).
\]
Choose $2-2/s<\alpha<2/q$ and
$r=2q(1-\alpha)/(2-q\alpha)>q$. Interpolation with
$\|\rho\dot u\|_{L^2}\le C A_2$ yields
\[
 \|\rho\dot u\|_{L^q}
 \le C_{q,\alpha,T}
       \bigl(A_2+A_2^\alpha\|\nabla\dot u\|_{L^2}^{1-\alpha}\bigr).
\]
The additive term lies in $L^s(0,S)$ because $s<2$. For the product term,
\[
\begin{aligned}
 &\int_0^S A_2^{s\alpha}\|\nabla\dot u\|_{L^2}^{s(1-\alpha)}\,\mathrm dt\\
 &\quad\le
 \left(\int_0^S A_2^2\,\mathrm dt\right)^{s\alpha/2}
 \left(\int_0^S t\|\nabla\dot u\|_{L^2}^2\,\mathrm dt
       \right)^{s(1-\alpha)/2}
 \left(\int_0^S t^{-s(1-\alpha)/(2-s)}\,\mathrm dt
       \right)^{1-s/2}\le C_{q,s,T},
\end{aligned}
\]
because $s(1-\alpha)/(2-s)<1$. Similarly,
\[
 \int_0^S tA_2^{2\alpha}\|\nabla\dot u\|_{L^2}^{2(1-\alpha)}\,\mathrm dt
 \le
 \left(\int_0^S tA_2^2\,\mathrm dt\right)^\alpha
 \left(\int_0^S t\|\nabla\dot u\|_{L^2}^2\,\mathrm dt
       \right)^{1-\alpha}\le C_T.
\]
On $\TT$,
$\int\rho\dot u\dd=0$ also gives
$\|\dot u\|_{H^1}\le C\|\nabla\dot u\|_{L^2}$.
\end{proof}

\begin{lemma}[Density gradients and strong elliptic bounds]\label{lem:gradient-closure}
Let $\Omega=\TT$ or $\R^2$, $q>2$, and $0<S\le T$.
Let $(\rho,u)$ be a local strong solution of \eqref{eq:system}, with
$\lambda$ satisfying \eqref{eq:lambda-regularity}, such that
\[
\begin{gathered}
 0\le\rho\le K,\qquad
 \sup_{t\le S}\bigl(\|\nabla u(t)\|_{L^2}+A_1(t)\bigr)
       +\int_0^S A_2^2\,dt\le K,\\
 \nabla\rho_0\in L^2\cap L^q,\qquad
 \rho\dot u\in L^1(0,S;L^q).
\end{gathered}
\]
Then
\begin{equation}\label{eq:density-gradient-bound}
 \sup_{t\le S}\|\nabla\rho(t)\|_{L^2\cap L^q}
       +\int_0^S\|\nabla u\|_{L^\infty}\,dt\le C_T,
\end{equation}
and, for almost every $t\in(0,S)$,
\begin{equation}\label{eq:elliptic-gradient-density}
\begin{aligned}
 \|\nabla^2u(t)\|_{L^2}&\le C_T(1+A_2(t)),\\
 \|\nabla^2u(t)\|_{L^q}&\le C_T\bigl(1+\|\rho\dot u(t)\|_{L^q}\bigr).
\end{aligned}
\end{equation}
The constant depends only on $\mu,\gamma,q,T,K$,
$\|\lambda\|_{C^1([0,K])}$, $\|\nabla\rho_0\|_{L^2\cap L^q}$, and
$\|\rho\dot u\|_{L^1(0,S;L^q)}$, and is independent of a positive
density lower bound.
\end{lemma}
\begin{proof}
Set $\Phi=(2\mu+\lambda(\rho))\nabla\rho$.
Differentiating continuity and the flux identity gives
\begin{equation}\label{eq:weighted-gradient-transport}
 \Dt\Phi+(\nabla u)^T\Phi+\Phi\diver u+\rho\nabla G
             +\frac{\gamma\rho^\gamma}{2\mu+\lambda(\rho)}\Phi=0.
\end{equation}
The $\rho\lambda'\nabla\rho\diver u$ terms cancel. Testing by
$|\Phi|^{q-2}\Phi$ gives
\[
 \frac{d}{dt}\|\Phi\|_{L^q}
 \le C\|\nabla u\|_{L^\infty}\|\Phi\|_{L^q}
                          +C\|\rho\dot u\|_{L^q}.
\]
Elliptic estimates, interpolation, and the flux definition yield
\[
\begin{aligned}
 \|\diver u\|_{L^\infty}+\|\omega\|_{L^\infty}
 &\le C\bigl(1+\|\rho\dot u\|_{L^q}^{q/(2(q-1))}\bigr),\\
 \|\nabla^2u\|_{L^q}
 &\le C\|\rho\dot u\|_{L^q}
       +C(1+\|\diver u\|_{L^\infty})\|\Phi\|_{L^q}.
\end{aligned}
\]
The Beale--Kato--Majda inequality
\cite[Lemma~2.4]{HuangLi2016} on $\TT$, and
\cite[Lemma~2.7]{HuangLi2022} on $\R^2$, imply
\[
 \|\nabla u\|_{L^\infty}
 \le C\bigl(1+\|\rho\dot u\|_{L^q}\bigr)
                         \log(e+\|\Phi\|_{L^q}).
\]
Since $q/(2(q-1))<1$, the logarithmic factor involving the force can be absorbed. Consequently,
\[
 \frac{d}{dt}\log\log(e+\|\Phi\|_{L^q})
                \le C\bigl(1+\|\rho\dot u\|_{L^q}\bigr).
\]
Integrating this inequality gives the $L^q$ gradient estimate, the Lipschitz bound,
and the $L^q$ elliptic estimate.

At exponent two, Hodge and interpolation give
\[
 \|\nabla^2u\|_{L^2}
 \le C_T\bigl(1+A_2+\|\nabla\rho\|_{L^2}\bigr)
              +\tfrac12\|\nabla^2u\|_{L^2}.
\]
Indeed, the coefficient-gradient term is bounded by
$C_T\|\diver u\|_{L^{2q/(q-2)}}$ and contributes a power $2/q<1$
of $\|\nabla^2u\|_{L^2}$. Differentiated continuity gives
\[
 \frac{d}{dt}\|\nabla\rho\|_{L^2}
 \le C\|\nabla u\|_{L^\infty}\|\nabla\rho\|_{L^2}
                                  +C\|\nabla^2u\|_{L^2}.
\]
Gronwall's inequality, together with $A_2\in L^2(0,S)$ and the Lipschitz bound,
completes the proof of \eqref{eq:density-gradient-bound} and the $L^2$ elliptic estimate.
\end{proof}

\subsection{Periodic strong estimates}
Throughout this subsection, $\Omega=\TT$.\label{per:sec:strong-estimates}

Fix $0<S\le T$ and $q>2$, with $\rho_0\in W^{1,q}$.
Constants $C_{q,T}$ may additionally depend on $\|\rho_0\|_{W^{1,q}}$
and $\sup_{[0,C]}|\lambda'|$, where $C$ is the density bound in
Proposition~\ref{per:prop:apriori}. No bound for $\lambda''$ at zero is used.
We also use $t\le\max\{1,T\}\min\{1,t\}$ on $[0,T]$.

\begin{proposition}[Finite-time strong estimates]\label{per:prop:strong-estimates}
Under the preceding assumptions, for every $1<s<q/(q-1)$,
\begin{equation}\label{per:eq:strong-estimates}
\begin{aligned}
 &\sup_{0\le t\le S}\bigl(\|\rho(t)\|_{W^{1,q}}
                              +\|\rho_t(t)\|_{L^{2}}\bigr)
 +\sup_{0<t\le S}t\|u(t)\|_{H^2}^2\\
 &\quad+\int_0^S\bigl(\|\nabla u\|_{L^{\infty}}
                 +\|u\|_{H^2}^2
                 +\|u\|_{W^{2,q}}^s
                 +\|\sqrt\rho\,u_t\|_{L^{2}}^2\bigr)\,\mathrm dt\\
 &\quad+\int_0^S t\bigl(\|u\|_{W^{2,q}}^2
                            +\|u_t\|_{H^1}^2\bigr)\,\mathrm dt
 \le C_{q,s,T}.
\end{aligned}
\end{equation}
\end{proposition}
\begin{proof}
Lemma~\ref{lem:gradient-closure} applies by
Proposition~\ref{per:prop:apriori} and Lemma~\ref{lem:force-integrability}.
Its Hessian bound and the force estimates give both $W^{2,q}$ time integrals.

The $L^2$ bound in \eqref{eq:elliptic-gradient-density}$_1$ and
the low-order $H^1$ bound give
\begin{equation}\label{per:eq:h2-velocity}
 \|u\|_{H^2}\le C_{q,T}(1+A_2),
\end{equation}
and Proposition~\ref{per:prop:apriori} gives the asserted time bounds.

For $u_t=\dot u-u\cdot\nabla u$, use
\[
\begin{aligned}
 \|u\cdot\nabla u\|_{L^{2}}^2
 &\le\|u\|_{L^{4}}^2\|\nabla u\|_{L^{4}}^2
 \le C_{q,T}(1+\|u\|_{H^2}^2),\\
 \|\nabla(u\cdot\nabla u)\|_{L^{2}}^2
 &\le C\|\nabla u\|_{L^{4}}^4
       +C\|u\|_{L^{2q/(q-2)}}^2\|\nabla^2u\|_{L^{q}}^2\\
 &\le C_{q,T}(1+A_2^2+\|\nabla^2u\|_{L^{q}}^2).
\end{aligned}
\]
Together with Lemma~\ref{lem:force-integrability}, these give the $u_t$ bounds.
Finally,
\[
 \|\rho_t\|_{L^{2}}
 \le\|u\|_{L^{2q/(q-2)}}\|\nabla\rho\|_{L^{q}}
                +\|\rho\|_{L^{\infty}}\|\diver u\|_{L^{2}}\le C_{q,T}.
\]
\end{proof}

\subsection{Periodic construction and uniqueness}
\label{per:sec:construction}

\begin{lemma}[Smooth coefficient approximations]\label{lem:coefficient-approximation}
Under \eqref{eq:lambda-regularity}, there are smooth positive coefficients
$\lambda_\varepsilon\to\lambda$ locally in $C^1$, with uniform
logarithmic derivative bounds and local $C^1$ bounds for $0<\varepsilon<1$.
If $\lambda$ satisfies \eqref{per:eq:primitive-growth}, the same
primitive growth is preserved with uniform constants.
\end{lemma}
\begin{proof}
Choose $0\le\chi\in C_c^\infty((-1,1))$ with integral one and set
\[
 \lambda_\varepsilon(r)=\varepsilon+
       \int_{-1}^1\chi(z)\lambda(r+\varepsilon(2+z))\,\mathrm dz.
\]
Then $\lambda_\varepsilon>0$, it converges locally in $C^1$, and
\[
 r|\lambda_\varepsilon'(r)|
       \le C_\lambda(2\mu+\lambda_\varepsilon(r)).
\]
For $r\ge3$, $r+\varepsilon(2+z)\in[r,2r]$, so
\eqref{eq:coefficient-bounds} gives
$2\mu+\lambda_\varepsilon(r)\ge2^{-b}(2\mu+\lambda(r))$. Hence
\[
 \int_1^r\frac{2\mu+\lambda_\varepsilon(z)}z\,\mathrm dz
       \ge c\theta(r)-C,\qquad r\ge3.
\]
The primitive growth follows by increasing its fixed threshold.
\end{proof}

\begin{lemma}[Global positive-density construction]\label{per:lem:positive-construction}
Smooth coefficients satisfying Assumption~\ref{per:ass:lambda} and smooth
periodic data with positive density generate a global positive strong
solution satisfying \eqref{per:eq:main-bound} and \eqref{per:eq:strong-estimates}.
\end{lemma}
\begin{proof}
For positive $H^2$ data, use the transport--parabolic iteration
\[
\begin{aligned}
 &\rho_t+\diver(\rho v)=0,\\
 &\rho u_t-\mu\Delta u
 -\nabla((\mu+\lambda(\rho))\diver u)\\
 &\qquad=-\rho v\cdot\nabla u-\nabla P,
\end{aligned}
\]
where $v$ is the preceding iterate. Transport preserves positivity;
composition and parabolic estimates give, on a common short interval,
\[
 \rho,u\in C([0,t_*];H^2),\qquad
 u\in L^2(0,t_*;H^3),\quad u_t\in L^2(0,t_*;H^1).
\]
The $L^2$ difference estimates give convergence after shortening the
interval. Its length depends on the $H^2$ data, the density lower bound,
and local coefficient bounds. This follows the framework of
Solonnikov~\cite{Solonnikov1980} and
\cite[Lemma~2.1 and Remark~2.1]{HuangLi2016}. For smooth general $\lambda$,
the same composition estimates on the positive density range suffice.

Suppose the maximal $H^2$ existence time $T^*$ is finite. Since
$H^2(\TT)\hookrightarrow W^{1,4}(\TT)$, apply
Propositions~\ref{per:prop:apriori} and \ref{per:prop:strong-estimates}
with $q=4$, uniformly for $S<T^*$.

Along the flow,
\begin{equation}\label{per:eq:positive-density-lower}
 \rho(t,x)\ge\inf\rho_0
 \exp\!\left(-\int_0^{T^*}\|\diver u\|_{L^{\infty}}\,\mathrm dt\right)>0.
\end{equation}
Here $C_*$ may depend on this lower bound and the higher initial norms.
Since $A_2(0)<\infty$, \eqref{per:eq:acceleration-differential} gives
\[
 \sup_{0\le t<S}A_2^2(t)
       +\int_0^S\|\nabla\dot u\|_{L^{2}}^2\,\mathrm dt\le C_*.
\]
Equation~\eqref{per:eq:h2-velocity}, weighted Poincar\'e, and elliptic
estimates give
\begin{equation}\label{per:eq:positive-flux-h2}
\begin{aligned}
 &\sup_{0\le t<S}\|u(t)\|_{H^2}\le C_*,\\
 &\int_0^S\bigl(\|\dot u\|_{H^1}^2
                  +\|\nabla^2G\|_{L^{2}}^2+\|\nabla^2\omega\|_{L^{2}}^2\bigr)
                  \,\mathrm dt\le C_*.
\end{aligned}
\end{equation}
Indeed,
\[
 \|\nabla(\rho\dot u)\|_{L^{2}}
 \le\|\nabla\rho\|_{L^{4}}\|\dot u\|_{L^{4}}
                    +\|\rho\|_{L^{\infty}}\|\nabla\dot u\|_{L^{2}}.
\]

For the density, the $W^{1,4}$ bound and
\[
 \nabla\Phi=(2\mu+\lambda)\nabla^2\rho
                   +\lambda'\nabla\rho\otimes\nabla\rho
\]
give
\begin{equation}\label{per:eq:positive-density-h2-comparison}
 \|\rho\|_{H^2}\le C_*(1+\|\nabla\Phi\|_{L^{2}}).
\end{equation}
To estimate $\nabla\Phi$, use the $L^4$ bound for $\nabla\rho$ and interpolation:
\[
\begin{aligned}
 \||\nabla\rho||\nabla G|\|_{L^{2}}
 &\le C_*(1+\|\nabla^2G\|_{L^{2}}^{1/2}),\\
 \||\nabla\rho||\nabla^2u|\|_{L^{2}}
 &\le C_*(1+\|\nabla^3u\|_{L^{2}}^{1/2}).
\end{aligned}
\]
The differentiated divergence relation gives
\[
 \|\nabla^3u\|_{L^{2}}
 \le C_*\bigl[\|\nabla^2G\|_{L^{2}}+\|\nabla^2\omega\|_{L^{2}}
       +(1+\|\nabla u\|_{L^{\infty}})(1+\|\nabla\Phi\|_{L^{2}})\bigr].
\]
Here second coefficient derivatives are used only on the positive density
range. Differentiate \eqref{eq:weighted-gradient-transport}, test in $L^2$,
and apply these product estimates to obtain
\[
 \frac{\mathrm d}{\mathrm dt}\|\nabla\Phi\|_{L^{2}}
 \le C_*(1+\|\nabla u\|_{L^{\infty}})(1+\|\nabla\Phi\|_{L^{2}})
             +C_*(\|\nabla^2G\|_{L^{2}}+\|\nabla^2\omega\|_{L^{2}}).
\]
Gronwall gives $\sup_{t<S}\|\rho(t)\|_{H^2}\le C_*$.
Since $\rho\dot u\in L^2L^4$, \eqref{eq:elliptic-gradient-density}$_2$
gives $\nabla u\in L^2L^\infty$ and hence $u\in L^2H^3$.
Also, $\|u\cdot\nabla u\|_{H^1}\le C\|u\|_{H^2}^2$ gives
$u_t\in L^2H^1$, uniformly for $S<T^*$.

In the transport equations for $\nabla\Phi$ and $\nabla\rho$, the coefficients
lie in $L^1L^\infty$ and the sources in $L^1L^2$ and $L^1L^4$, respectively.
Thus $\nabla\Phi\in C_tL^2$, $\rho\in C_tW^{1,4}$, and the identity
for $\nabla\Phi$ gives $\rho\in C([0,T^*];H^2)$.
The Hilbert triple $H^3\subset H^2\subset H^1$ gives the same continuity
for $u$. Restarting from the positive $H^2$ terminal data contradicts maximality.
\end{proof}

\begin{proof}[Existence in Theorem~\ref{per:thm:main}]
Let $j_\varepsilon$ be a standard periodic mollifier and set
\[
 \rho_{0,\varepsilon}
 =\frac{j_\varepsilon*\rho_0+\varepsilon}{1+\varepsilon},
 \qquad u_{0,\varepsilon}=j_\varepsilon*u_0.
\]
These positive, unit-mass data converge in $W^{1,q}\times H^1$.
Lemmas~\ref{lem:coefficient-approximation} and
\ref{per:lem:positive-construction} give global solutions
$(\rho_\varepsilon,u_\varepsilon)$ satisfying \eqref{per:eq:main-bound},
\eqref{eq:force-integrability}, and \eqref{per:eq:strong-estimates}
uniformly on $[0,T]$. The local $\lambda_\varepsilon'$ bounds are also uniform.

The bounds for $\rho_\varepsilon$ in $L^\infty W^{1,q}$ and for its time
derivative in $L^\infty L^2$ give, along a subsequence,
\begin{equation}\label{per:eq:density-compactness}
 \rho_\varepsilon\longrightarrow\rho
 \quad\text{in }C([0,T];C^\alpha(\TT)),
 \qquad 0<\alpha<1-2/q.
\end{equation}
The velocity bounds in $L^2(\tau,T;W^{2,q})$ and time-derivative bounds in
$L^2(\tau,T;H^1)$ give strong convergence in $L^2(\tau,T;H^1)$. The $L^\infty H^1$
bound controls the omitted initial interval, giving
\begin{equation}\label{per:eq:velocity-compactness}
 u_\varepsilon\longrightarrow u
 \quad\text{in }L^2(0,T;H^1).
\end{equation}
Moreover, $\rho_\varepsilon u_\varepsilon$ is bounded in $L^\infty H^1$,
and the momentum equation bounds its time derivative in $L^\infty H^{-1}$, since
the stress and convection tensor are bounded in $L^\infty L^2$. Thus
\begin{equation}\label{per:eq:momentum-compactness}
 \rho_\varepsilon u_\varepsilon\longrightarrow\rho u
 \quad\text{in }C([0,T];L^2).
\end{equation}
Together with uniform convergence of $\lambda_\varepsilon(\rho_\varepsilon)$,
these convergences identify all terms in \eqref{eq:system}.

The differentiated continuity equation has $L^1L^\infty$ coefficients
and source $-\rho\nabla\diver u\in L^1L^q$. Its transport formula gives
$\rho\in C([0,T];W^{1,q})$ with initial value $\rho_0$.
Equation~\eqref{per:eq:momentum-compactness} gives the momentum trace.
For $\delta>0$, $\rho u/\sqrt{\rho+\delta}\in C([0,T];L^2)$, and
\[
 \left\|\frac{\rho u}{\sqrt{\rho+\delta}}
                        -\sqrt\rho\,u\right\|_{L^{2}}
 \le\sqrt\delta\,\|u\|_{L^{2}}.
\]
Letting $\delta\downarrow0$ uniformly proves continuity of
$\sqrt\rho\,u$ and identifies its trace as $\sqrt{\rho_0}u_0$.

On $[\tau,T]$, the $L^2H^1$ bounds and strong velocity convergence identify
the limits of $\partial_tu_\varepsilon$ and $\dot u_\varepsilon$.
Uniform density convergence identifies their density-weighted weak limits,
including the time-weighted acceleration. Weak lower semicontinuity, followed by
$\tau\downarrow0$, gives \eqref{per:eq:main-bound} and
\eqref{per:eq:strong-estimates}. Also $u_t\in L^2(\tau,T;H^1)$ gives
$u\in C([\tau,T];H^1)$, completing \eqref{per:eq:strong-class}.
A diagonal extraction over integer time horizons and exponents
$s\uparrow q/(q-1)$ gives one global solution with all the stated
bounds.
\end{proof}

\begin{proof}[Uniqueness in Theorem~\ref{per:thm:main}]
Let $(\rho_1,u_1)$ and $(\rho_2,u_2)$ have the same initial density
and momentum. Constants may depend on their finite norms on $[0,T]$.
On the common density range, $P$ and $\lambda$ are Lipschitz.
Weighted Poincar\'e and Sobolev give
\[
\begin{aligned}
 &\|u_1-u_2\|_{H^1}+\|u_1-u_2\|_{L^{2q/(q-2)}}\\
 &\qquad\le C_{q,T}\bigl(
   \|\sqrt{\rho_1}(u_1-u_2)\|_{L^{2}}+\|\nabla(u_1-u_2)\|_{L^{2}}\bigr).
\end{aligned}
\]
The density difference satisfies
\[
\begin{aligned}
 &\partial_t(\rho_1-\rho_2)+u_2\cdot\nabla(\rho_1-\rho_2)
                       +(\rho_1-\rho_2)\diver u_2\\
 &\qquad=-(u_1-u_2)\cdot\nabla\rho_1-\rho_1\diver(u_1-u_2).
\end{aligned}
\]
Since $\nabla u_2\in L^1L^\infty$ and the initial density difference is zero,
\begin{equation}\label{per:eq:difference-density}
\begin{aligned}
 &\|\rho_1(t)-\rho_2(t)\|_{L^{2}}^2\\
 &\quad\le C_{q,T}t\int_0^t\bigl(
  \|\sqrt{\rho_1}(u_1-u_2)\|_{L^{2}}^2
                +\|\nabla(u_1-u_2)\|_{L^{2}}^2\bigr)\,\mathrm ds.
\end{aligned}
\end{equation}
The velocity difference equation is
\[
\begin{aligned}
 &\rho_1\bigl(\partial_t(u_1-u_2)
                 +u_1\cdot\nabla(u_1-u_2)\bigr)-\mu\Delta(u_1-u_2)\\
 &\quad-\nabla\bigl((\mu+\lambda(\rho_1))\diver(u_1-u_2)\bigr)\\
 &\qquad=-(\rho_1-\rho_2)\dot u_2
                    -\rho_1(u_1-u_2)\cdot\nabla u_2\\
 &\qquad\quad+\nabla\bigl((\lambda(\rho_1)-\lambda(\rho_2))\diver u_2\bigr)
                    -\nabla(P(\rho_1)-P(\rho_2)).
\end{aligned}
\]
Testing by $u_1-u_2$, using weighted Poincar\'e and the Lipschitz bounds
for $P$ and $\lambda$, gives
\begin{equation}\label{per:eq:difference-energy}
\begin{aligned}
 &\frac{\mathrm d}{\mathrm dt}
       \|\sqrt{\rho_1}(u_1-u_2)\|_{L^{2}}^2
                     +c\|\nabla(u_1-u_2)\|_{L^{2}}^2\\
 &\quad\le C_{q,T}(1+\|\nabla u_2\|_{L^{\infty}})
                    \|\sqrt{\rho_1}(u_1-u_2)\|_{L^{2}}^2\\
 &\qquad+C_{q,T}\bigl(1+\|\nabla u_2\|_{L^{\infty}}^2
                     +\|\dot u_2\|_{L^{q}}^2\bigr)\|\rho_1-\rho_2\|_{L^{2}}^2.
\end{aligned}
\end{equation}
By \eqref{per:eq:strong-class}, Sobolev, and
$\|u_2\cdot\nabla u_2\|_{L^q}
\le C_q\|u_2\|_{H^1}\|u_2\|_{W^{2,q}}$,
\[
 \int_0^T\left[1+\|\nabla u_2\|_{L^{\infty}}
       +t\bigl(\|\nabla u_2\|_{L^{\infty}}^2+\|\dot u_2\|_{L^{q}}^2\bigr)\right]
       \,\mathrm dt<\infty.
\]
The initial weighted-velocity traces imply
$\|\sqrt{\rho_1}(u_1-u_2)\|_{L^{2}}\to0$ as $t\downarrow0$, since
\[
 \|(\sqrt{\rho_1}-\sqrt{\rho_2})u_2\|_{L^{2}}^2
 \le\|\rho_1-\rho_2\|_{L^{\infty}}\|u_2\|_{L^{2}}^2\longrightarrow0.
\]
Set
\[
\begin{aligned}
 \mathcal Y(t)={}&\sup_{0<s\le t}
    \|\sqrt{\rho_1(s)}(u_1(s)-u_2(s))\|_{L^{2}}^2\\
 &+\int_0^t\|\nabla(u_1-u_2)\|_{L^{2}}^2\,\mathrm ds.
\end{aligned}
\]
By \eqref{per:eq:difference-density},
$\|\rho_1(t)-\rho_2(t)\|_{L^{2}}^2\le C_{q,T}t\mathcal Y(t)$.
The right-hand side of \eqref{per:eq:difference-energy} is therefore integrable.
Integrating from a positive time and then letting that time tend to zero gives
\[
 \mathcal Y(t)\le C_{q,T}\int_0^t
 \left[1+\|\nabla u_2\|_{L^{\infty}}
     +s\bigl(1+\|\nabla u_2\|_{L^{\infty}}^2+\|\dot u_2\|_{L^{q}}^2\bigr)\right]
 \mathcal Y(s)\,\mathrm ds.
\]
Gronwall gives $\rho_1=\rho_2$ and $u_1=u_2$.
\end{proof}

\subsection{Whole-space construction and continuation}
Throughout this subsection, $\Omega=\R^2$.\label{ws:sec:construction}

The local construction uses \eqref{eq:lambda-regularity} alone.

\begin{lemma}[Local construction]
\label{ws:lem:general-local}
Assume \eqref{eq:lambda-regularity}.
Data satisfying \eqref{ws:eq:initial-data} admit a unique strong solution
on some $[0,T_0]$. It satisfies \eqref{ws:eq:strong-class}, with $T=T_0$, and
$\nabla u\in L^1(0,T_0;L^\infty)$.
The time $T_0$ depends only on the initial norms, a ball containing a
fixed positive fraction of the initial mass, $\mu,\gamma,a,q,\eta_0$,
$C_\lambda$, and $\|\lambda\|_{C^1([0,1])}$.
It is independent of the density lower bound and of $\lambda''$ near vacuum.
The time-integrated and time-weighted bounds in \eqref{ws:eq:strong-class}
have the same data dependence.
No initial acceleration compatibility condition is required.
The solution also has the initial upper trace
\begin{equation}\label{ws:eq:initial-gradient-energy}
\begin{aligned}
 &\mathop{\mathrm{ess\,limsup}}_{t\downarrow0}
       \int(\mu(\curl u)^2+(2\mu+\lambda(\rho)) (\diver u)^2-2P\diver u)\dd\\
 &\qquad\le\int\bigl(
       \mu|\curl u_0|^2
       +(2\mu+\lambda(\rho_0))|\diver u_0|^2
       -2\rho_0^\gamma\diver u_0\bigr)\dd .
\end{aligned}
\end{equation}
This is an energy upper trace, not an assertion that
$\nabla u(t)\to\nabla u_0$ strongly at vacuum.
\end{lemma}
\begin{proof}
We adapt \cite[Sections~2--3 and~5]{LiLiang2014} from power laws to
general coefficients. The logarithmic derivative bound gives, for some $b\ge1$,
\[
 \lambda(r)\le C(1+r^b),\qquad
 \sup_{0\le r\le K}|\lambda'(r)|\le C(1+K^{b-1})\quad(K\ge1).
\]
Include $\lambda(0)$ in the constant-coefficient operator and use
$\lambda(\rho)-\lambda(0)$ in all integrable coefficient norms.
For $0\le\rho\le K$ and $L_K=\sup_{[0,K]}|\lambda'|$,
\[
 |\lambda(\rho)-\lambda(0)|\le L_K\rho,\qquad
 |\rho\lambda'(\rho)|\le L_K\rho,\qquad
 |\nabla\lambda(\rho)|\le L_K|\nabla\rho|.
\]
The coercive coefficient remains at least $2\mu$.

The transport identity replacing (3.14) of \cite{LiLiang2014} is
\[
\begin{aligned}
 \lambda_t
 &=-\diver\bigl((\lambda(\rho)-\lambda(0))u\bigr)
          +(\lambda(\rho)-\lambda(0)-\rho\lambda'(\rho))\diver u\\
 &=-u\cdot\nabla\lambda-\rho\lambda'(\rho)\diver u.
\end{aligned}
\]
On the approximating balls, $u\cdot n=0$ and integration by parts gives
\[
\begin{aligned}
 \int\lambda_t(\diver u)^2\dd
 &=2\int(\lambda(\rho)-\lambda(0))\,
                         u\cdot\nabla(\diver u)\,\diver u\dd\\
 &\quad+\int(\lambda(\rho)-\lambda(0)
                   -\rho\lambda'(\rho))(\diver u)^3\dd.
\end{aligned}
\]
These terms satisfy the density-weighted estimate (3.13) there.
In the time-differentiated equation,
\[
 \int|\rho\lambda'|(\diver u)^2|\diver u_t|\dd
 \le\varepsilon\mu\|\nabla u_t\|_{L^2}^2
      +C_\varepsilon\|\rho\lambda'\|_{L^\infty}^2\|\diver u\|_{L^4}^4.
\]
The weighted H\"older estimate (3.30) in that paper controls $u\cdot\nabla\lambda$, since
$\|\bar x^{a/2}\nabla\lambda\|_{L^2\cap L^q}
 \le L_K\|\bar x^{a/2}\nabla\rho\|_{L^2\cap L^q}$.
The weighted density estimate (3.37) and the chain rule give the
weighted gradient bounds for $\lambda$, replacing (3.38).
Thus the proofs of (3.22)--(3.24) and (3.39) require only first coefficient
derivatives, with bounds polynomial in the density norms.

Use the positive coefficients of Lemma~\ref{lem:coefficient-approximation};
the preceding bounds are uniform for
$\lambda_\varepsilon(\rho)-\lambda_\varepsilon(0)$.
Let $j_\delta$ be a nonnegative standard mollifier and let
$0\le\chi_R\le1$ equal one on $B_{R/4}$ and vanish outside $B_{R/2}$,
with $|\nabla^k\chi_R|\le C_kR^{-k}$.
For large $R$, put
\[
 \rho_{0,R}=(\chi_R\rho_0)*j_{1/R}
                    +R^{-1}e^{-|x|^2},\qquad
 \varepsilon=R^{-1}.
\]
The densities are smooth, strictly positive on $\overline B_R$, and
$\bar x^a\rho_{0,R}\to\bar x^a\rho_0$ in
$L^1\cap H^1\cap W^{1,q}$ on $\R^2$.
They retain a fixed positive mass in one fixed ball.
The compact $C^1$ convergence of the coefficients also gives
\[
 \begin{aligned}
 (\lambda_\varepsilon(\rho_{0,R})-\lambda_\varepsilon(0))&\longrightarrow (\lambda(\rho_0)-\lambda(0))
                                      &&\text{in }L^2,\\
 \bar x^{a/2}\nabla\lambda_\varepsilon(\rho_{0,R})
      &\longrightarrow\bar x^{a/2}\nabla\lambda(\rho_0)
                                      &&\text{in }L^2\cap L^q.
 \end{aligned}
\]
The first limit follows from
$\lambda_\varepsilon(r)-\lambda_\varepsilon(0)=\int_0^r\lambda_\varepsilon'(s)\,ds$.

Choose vector fields $v_{i,R}\in C_c^\infty(B_R;\R^2)$ with
$v_{i,R}\to\partial_i u_0$ in $L^2(\R^2)$, and set
$f_R=(\sqrt{\rho_0}u_0)*j_{1/R}$.
Define $u_{0,R}$ by the coefficient-independent elliptic problem
\[
 \begin{cases}
 -\Delta u_{0,R}+(\rho_{0,R}+R^{-1})u_{0,R}
       =\sqrt{\rho_{0,R}}f_R-\displaystyle\sum_{i=1}^2
                                       \partial_i v_{i,R},&x\in B_R,\\
 u_{0,R}\cdot n=0,\quad\curl u_{0,R}=0,&x\in\partial B_R.
 \end{cases}
\]
This is the approximation in \cite[(5.5)]{LiLiang2014}.
Testing it and applying the slip Hodge estimate give
\[
 R^{-1}\|u_{0,R}\|_{L^2}^2+\|\sqrt{\rho_{0,R}}u_{0,R}\|_{L^2}^2
                      +\|\nabla u_{0,R}\|_{L^2}^2\le C.
\]
If $\varphi_R$ equals one on $B_{R/2}$ and vanishes near
$\partial B_R$, with $|\nabla^k\varphi_R|\le C_kR^{-k}$, then
\[
 \nabla(\varphi_Ru_{0,R})\longrightarrow\nabla u_0,\qquad
 \sqrt{\rho_{0,R}}\,\varphi_Ru_{0,R}
                    \longrightarrow\sqrt{\rho_0}u_0
                    \quad\text{in }L^2(\R^2).
\]
For any weak elliptic limit $w$, testing its difference equation by
bounded velocity truncations and logarithmic spatial cutoffs gives
\[
 -\Delta(w-u_0)+\rho_0(w-u_0)=0
 \quad\Longrightarrow\quad
 \|\nabla(w-u_0)\|_{L^2}^2+\|\sqrt{\rho_0}(w-u_0)\|_{L^2}^2=0.
\]
Thus every weak limit is $u_0$, since $M>0$. Moreover,
\[
\begin{aligned}
 &\int_{B_R}\bigl(
       |\diver u_{0,R}|^2+|\curl u_{0,R}|^2
             +(\rho_{0,R}+R^{-1})|u_{0,R}|^2\bigr)\dd\\
 &=\int_{B_R}\bigl(
       \sqrt{\rho_{0,R}}f_R\cdot u_{0,R}
                    +\sum_i v_{i,R}\cdot\partial_i u_{0,R}\bigr)\dd\\
 &\longrightarrow
       \|\nabla u_0\|_{L^2}^2+\|\sqrt{\rho_0}u_0\|_{L^2}^2 .
\end{aligned}
\]
Lower semicontinuity gives strong convergence of the zero-extended
divergences, curls, and weighted velocities, and
$R^{-1}\|u_{0,R}\|_{L^2}^2\to0$.
The cutoff gradient error is bounded by
$CR^{-1}\|u_{0,R}\|_{L^2}\le CR^{-1/2}$; the whole-space Hodge identity
then gives the asserted gradient convergence.

Solve the smooth positive-density problem on $B_R$, with coefficient
$\lambda_{1/R}$, the same slip conditions, and damping $-R^{-1}u_R$.

Let $\mathcal Q_R(t)$ be one plus the instantaneous norms
$\|\sqrt\rho u\|_{L^2}+\|\nabla u\|_{L^2}+R^{-1/2}\|u\|_{L^2}
+\|\bar x^a\rho\|_{L^1\cap H^1\cap W^{1,q}}$.
For a fixed $\alpha>1$, the composition bounds give
\[
 \|\lambda(\rho)-\lambda(0)\|_{L^2}
 +\|\bar x^{a/2}\nabla\lambda(\rho)\|_{L^2\cap L^q}
 +\|\rho\lambda'(\rho)\|_{L^\infty}
 \le C\mathcal Q_R^\alpha.
\]
Uniformly for $\lambda_\varepsilon$, the preceding adaptations of
\cite[Proposition~3.1]{LiLiang2014} yield
\[
 \mathcal Q_R(t)\le
 C\exp\!\left\{C\exp\!\left(C\int_0^t\mathcal Q_R(s)^\alpha\,ds
                                     \right)\right\},\qquad\alpha>1,
\]
together with the time-integrated and time-weighted velocity bounds
in \eqref{ws:eq:strong-class}.
The constants are independent of $R$ and the minimum density.
Choose $K>2C\exp(C\exp C)$ and $T_0\le K^{-\alpha}$ below the
uniform mass-localization time. Before the first time at which $\mathcal Q_R=K$,
$\int_0^{T_0}\mathcal Q_R^\alpha\,dt\le1$, so the displayed inequality gives
$\mathcal Q_R<K/2$. Hence $T_0$ is common to all approximations, and
$\int_0^{T_0}\|\nabla u_R\|_{L^\infty}\,dt\le C$. In particular,
\[
 \inf_{B_R}\rho_R(t)\ge
 \inf_{B_R}\rho_{0,R}
       \exp\!\left(-\int_0^{T_0}\|\diver u_R\|_{L^\infty}\,dt\right)>0.
\]
Classical higher-order continuation therefore extends each approximation
to $T_0$; its higher-order constants may depend on $R$.

Extend $\widetilde\rho_R=\varphi_R^4\rho_R$ and $w_R=\varphi_Ru_R$
by zero outside $B_R$.
The weighted density bounds control the cutoff errors, while
$|\lambda_\varepsilon(r)-\lambda_\varepsilon(0)|\le Cr$ and
$|\lambda_\varepsilon'(r)|\le C$ on the bounded density range control
the extended viscosity.
As in \cite[(5.11)--(5.17)]{LiLiang2014}, continuity, the velocity
bounds, and localized mass give, on each fixed ball,
\[
 \partial_t\widetilde\rho_R\in L^2(0,T_0;L^q),\qquad
 \partial_t w_R\in L^2(\tau,T_0;L^2),\quad \tau>0,
\]
uniformly in $R$. Local compactness yields
\[
 \widetilde\rho_R\longrightarrow\rho
       \quad\text{in }C([0,T_0];C_{\rm loc}(\R^2)),\qquad
 w_R\longrightarrow u
       \quad\text{in }L^2(\tau,T_0;H^1_{\rm loc}(\R^2)).
\]
Compact $C^1$ convergence identifies both
$\lambda_{1/R}(\rho_R)\diver u_R$ and
$\lambda'_{1/R}(\rho_R)\nabla\rho_R$ in the weak limit.
Convection converges on $(\tau,T_0)$ and has integral at most $C\tau$
on $(0,\tau)$; the damping vanishes in distributions. Compactly supported
tests therefore give the equations and initial momentum $\rho_0u_0$.

To obtain the initial upper trace, test momentum by $u_{R,t}$,
integrate the pressure term in time, and use the coefficient bounds and
\cite[(3.18), (3.23)]{LiLiang2014}:
\[
\begin{aligned}
 \frac{d}{dt}\biggl\{
   &\int_{B_R}(\mu(\curl u_R)^2+(2\mu+\lambda_{1/R}(\rho_R))(\diver u_R)^2-2P_R\diver u_R)\dd
                                  +R^{-1}\|u_R\|_{L^2}^2\biggr\}\\
 &\quad+c\|\sqrt{\rho_R}u_{R,t}\|_{L^2}^2\le C_{\rm loc}.
\end{aligned}
\]
Here $P_R=\rho_R^\gamma$, and $C_{\rm loc}$ depends on the local initial
norms but not on $R$ or the density lower bound. The initial integral
converges by the strong div--curl convergence and uniform initial
density convergence. Local density convergence and weighted mass give
global strong $L^2$ convergence of the pressure, since
\[
 \sup_{t\le T_0}
   \|P_R-\widetilde\rho_R^\gamma\|_{L^2(B_R)}^2
 \le C\sup_{t\le T_0}\int_{B_R\setminus B_{R/2}}\rho_R\dd
 \le C R^{-a}\longrightarrow0.
\]
After time integration, the pressure cross term converges and the
positive quadratic terms pass by lower semicontinuity. Letting
$t\downarrow0$ gives \eqref{ws:eq:initial-gradient-energy}, since
$C_{\rm loc}t\to0$.

For the kinetic trace, H\"older's inequality and the fixed-mass Hardy estimate give
\begin{equation}\label{ws:eq:local-kinetic-tail}
\begin{aligned}
 &\int\widetilde\rho_R\bar x^{(a-1)/2}|w_R|^2\dd\\
 &\quad\le
 \left(\int\widetilde\rho_R\bar x^a\dd\right)^{(a+1)/(2a)}
 \|\widetilde\rho_R\|_{L^\infty}^{(a-1)/(2a)}
 \|w_R\bar x^{-1/2}\|_{L^{4a/(a-1)}}^2\le C .
\end{aligned}
\end{equation}
The Hardy condition is strict because $4a/(a-1)>4$.
The limit has the same uniform tail bound. Local $H^1$ compactness,
density continuity, and the initial momentum trace imply
$\sqrt\rho u(t)\to\sqrt{\rho_0}u_0$ in $L^2$ as $t\downarrow0$.
Distributional momentum continuity allows the same argument at every
reference time, giving both momentum continuities in
\eqref{ws:eq:strong-class}. The Lipschitz transport equation, with
$\nabla^2u\in L^1_t(L^2\cap L^q)$, gives
$\rho\in C_t(H^1\cap W^{1,q})$; weighted mass gives $C_tL^1$ continuity.
The other bounds pass by lower semicontinuity, and uniqueness follows
from Subsection~\ref{ws:sec:uniqueness}.
\end{proof}

We next prove the bounds needed to restart the local solution.
Here $C_T$ depends on the fixed finite time, initial weighted norms,
parameters, coefficient, and \eqref{ws:eq:uniform-core}, but not on
the interval endpoint or density lower bound. We justify the calculations
on the whole-space solution using spatial cutoffs and time
difference quotients on $[\delta,t]$, then let $\delta\downarrow0$.

\begin{lemma}[Finite-time strong estimates]\label{ws:lem:finite-strong}
Let a local strong solution of \eqref{eq:system} satisfy the basic energy
bound and the estimates in \eqref{ws:eq:uniform-core}, with all
suprema and integrals restricted to its interval of existence and
with constants independent of its right endpoint.
For every finite $T$, the following bounds hold up to the lifespan endpoint:
\begin{equation}\label{ws:eq:finite-strong-bound}
\begin{aligned}
 &\sup_{0<t\leq T}\left\{
   \|\bar x^a\rho(t)\|_{L^1\cap H^1\cap W^{1,q}}
   +t\|\nabla^2u(t)\|_{L^2}^2
   +t\|\sqrt\rho u_t(t)\|_{L^2}^2\right\}\\
 &\quad+\int_0^T\Bigl(
   \|\nabla u\|_{L^\infty}+\|\nabla^2u\|_{L^2}^2
   +\|\nabla^2u\|_{L^q}^{(q+1)/q}
   +t\|\nabla^2u\|_{L^q}^2\Bigr)\,\mathrm dt\\
 &\quad+\int_0^T\Bigl(\|\sqrt\rho u_t\|_{L^2}^2
   +t\|\nabla u_t\|_{L^2}^2
   +t\|\bar x^{-1}u_t\|_{L^2}^2\Bigr)\,\mathrm dt
 \leq C_T.
\end{aligned}
\end{equation}
Here and below the upper endpoint is restricted to the solution's
interval of existence if that interval ends before $T$.
The mass lower bound \eqref{ws:eq:finite-mass-localization}
also holds on the same interval.
\end{lemma}

\begin{proof}
Set $p=(q+1)/q$ and apply Lemma~\ref{lem:force-integrability} with $s=p$.
We also use the weighted Hardy inequalities
\cite[Lemma~2.4, (2.7)--(2.8)]{LiLiang2014} with weights $\bar x^{-1}$ in
$L^2$ and $\bar x^{-\varepsilon}$ in $L^r$, where
$0<\varepsilon<1$ and $r>2/\varepsilon$.
Apply Lemma~\ref{lem:gradient-closure}, using
\eqref{ws:eq:uniform-core} and \eqref{eq:force-integrability}.
The pointwise elliptic bound in that lemma also yields
\begin{equation}\label{ws:eq:construction-u-q}
 \int_0^T\left(\|\nabla^2u\|_{L^q}^p
              +t\|\nabla^2u\|_{L^q}^2\right)\,dt\le C_T.
\end{equation}

The $L^2$ conclusions of Lemma~\ref{lem:gradient-closure} and
\eqref{ws:eq:uniform-core} give
\begin{equation}\label{ws:eq:construction-u-two}
 \sup_{t\leq T}\left(\|\nabla\rho(t)\|_{L^2}
                         +t\|\nabla^2u(t)\|_{L^2}^2\right)
             +\int_0^T\|\nabla^2u\|_{L^2}^2\,\mathrm dt\leq C_T.
\end{equation}

To propagate the weighted density, Hardy and Sobolev interpolation give
\begin{equation}\label{ws:eq:construction-weighted-velocity}
 \|u\bar x^{-\varepsilon}\|_{L^\infty}
 \leq C_{T,\varepsilon}(1+\|\nabla u\|_{L^q}),
 \qquad 0<\varepsilon<1.
\end{equation}
Indeed, differentiate $u\bar x^{-\varepsilon}$ and interpolate with
$\|u\bar x^{-\varepsilon}\|_{L^{4/\varepsilon}}$.
The right-hand side is integrable by \eqref{ws:eq:construction-u-two}.
Also,
\begin{equation}\label{ws:eq:construction-weighted-transport}
 \Dt(\bar x^a\rho)=\bigl(a u\cdot\nabla\log\bar x-\diver u\bigr)\bar x^a\rho.
\end{equation}
Its coefficient belongs to $L^1_tL^\infty_x$, so
\[
 \sup_{t\leq T}\|\bar x^a\rho(t)\|_{L^\infty}\leq C_T .
\]
Differentiate \eqref{ws:eq:construction-weighted-transport} and use
$|\nabla\log\bar x|\leq C\langle x\rangle^{-1}$,
$|\nabla^2\log\bar x|\leq C\langle x\rangle^{-2}$.
For $r=2,q$ this gives
\[
\begin{aligned}
 \frac{\mathrm d}{\mathrm dt}\|\nabla(\bar x^a\rho)\|_{L^r}
 \leq{}&C\bigl(\|\nabla u\|_{L^\infty}
                    +\|u\bar x^{-1/4}\|_{L^\infty}\bigr)\|\nabla(\bar x^a\rho)\|_{L^r}\\
       &+C_T\bigl(\|\nabla^2u\|_{L^r}+\|\nabla u\|_{L^r}
                            +\|u\bar x^{-1/4}\|_{L^\infty}\bigr).
\end{aligned}
\]
All coefficients and sources are integrable. Gronwall and
\eqref{ws:eq:finite-weighted-mass} yield
\begin{equation}\label{ws:eq:construction-weighted-density}
 \sup_{t\leq T}\|\bar x^a\rho(t)\|_{L^1\cap H^1\cap W^{1,q}}
 \leq C_T.
\end{equation}

The Eulerian time-derivative bounds follow from uniform local restarts.
By \eqref{ws:eq:construction-weighted-density},
\eqref{ws:eq:finite-mass-localization}, and the energy bound, the data at all
admissible restart times on $[0,S]$, $S\le T$, are uniformly bounded in
the norms required by Lemma~\ref{ws:lem:general-local}. Its quantitative estimates give
a common $0<\delta_T\le\min\{1,T_0\}$ and, with $\sigma=t-t_0$,
\[
\begin{aligned}
 &\sup_{0<\sigma\le\delta_T}
       \sigma\|\sqrt\rho u_t(t_0+\sigma)\|_{L^2}^2\\
 &\quad+\int_0^{\delta_T}\Bigl(
       \|\sqrt\rho u_t(t_0+\sigma)\|_{L^2}^2
       +\sigma\|\nabla u_t(t_0+\sigma)\|_{L^2}^2
       +\sigma\|\bar x^{-1}u_t(t_0+\sigma)\|_{L^2}^2\Bigr)
       \,\mathrm d\sigma\le C_T.
\end{aligned}
\]
The constants depend only on those common data bounds.
Uniqueness in \eqref{ws:eq:strong-class} identifies each restarted solution
with the original one on their overlap. This uses only finite local class norms.

Choose admissible $t_j$ within $\delta_T/16$ of grid points spaced
$\delta_T/4$ apart. At most $C(1+T/\delta_T)$ interior windows
\[
 I_j=[t_j+\delta_T/4,t_j+3\delta_T/4]\cap[0,S]
\]
cover $[\delta_T/2,S]$. On each window, $\sigma\ge\delta_T/4$ and
$t/\sigma\le4T/\delta_T$. Thus the displayed local estimates give
the unweighted $\|\sqrt\rho u_t\|_{L^2}^2$ integral and the required
supremum and integrals weighted by absolute time on $I_j$.
On $[0,\min\{\delta_T/2,S\}]$, use the original local estimates with
$\sigma=t$. Summing over the finite cover completes
\eqref{ws:eq:finite-strong-bound}, uniformly in the lifespan endpoint.
\end{proof}

\begin{proposition}[Continuation]\label{ws:prop:strong-continuation}
Under the structural condition~\ref{per:ass:lambda}, the local strong
solution with data \eqref{ws:eq:initial-data} extends to every finite time
and has the regularity \eqref{ws:eq:strong-class}.
\end{proposition}

\begin{proof}
If the maximal time $T_*$ were finite,
Lemma~\ref{ws:lem:finite-strong} and \eqref{ws:eq:uniform-core} would give
\[
 \sup_{t<T_*}\left(
 \|\bar x^a\rho(t)\|_{L^1\cap H^1\cap W^{1,q}}
       +\|\nabla u(t)\|_{L^2}+\|\sqrt\rho u(t)\|_{L^2}\right)<\infty,
\]
and the mass lower bound in
\eqref{ws:eq:finite-mass-localization} holds on one fixed ball.

Lemma~\ref{ws:lem:general-local} gives a common lifespan
for restarts at admissible $t_0\in[T_*/2,T_*)$.
For $t_0$ sufficiently close to $T_*$, the restarted solution extends
beyond $T_*$ and agrees with the original one by uniqueness.
This contradicts maximality, without imposing terminal compatibility
at vacuum. Hence $T_*=\infty$, and the local construction and
Lemma~\ref{ws:lem:finite-strong} give \eqref{ws:eq:strong-class}.
\end{proof}

\begin{remark}\label{ws:rem:construction-scope}
The density-gradient and weighted bounds may grow with the finite
time $T$; \eqref{ws:eq:uniform-core} is uniform in time.
No bound for $\lambda''$ near vacuum is asserted, nor is classical regularity.
\end{remark}

\subsection{Whole-space uniqueness}\label{ws:sec:uniqueness}

Let $(\rho,u)$ and $(\widetilde\rho,\widetilde u)$ satisfy
\eqref{ws:eq:strong-class} on $[0,T]$ with the same initial density
and momentum. Their positive mass and weighted density bounds give
uniform mass localization. Thus \cite[Lemma~2.4]{LiLiang2014} yields
\begin{equation}\label{ws:eq:uniqueness-hardy}
 \|v\bar x^{-\eta}\|_{L^p}
 \le C_{T,\eta,p}
       \bigl(\|\sqrt\rho\,v\|_{L^2}+\|\nabla v\|_{L^2}\bigr),
 \qquad 0<\eta<1,\quad p>\frac2\eta,
\end{equation}
and the same bound with $\widetilde\rho$ in place of $\rho$.
Constants may depend on $T$, but not on a density lower bound.

Put
\[
 H=\rho-\widetilde\rho,\qquad U=u-\widetilde u,\qquad
 r=\frac{1+a}{2},\qquad e(t)=\int\rho|U|^2\,\dd.
\]
These quantities are finite; indeed, \eqref{ws:eq:uniqueness-hardy} gives
\[
 \int |H||\widetilde u|^2\,\dd
 \le \|\bar x^rH\|_{L^2}
                         \|\widetilde u\bar x^{-r/2}\|_{L^4}^2
 \le C_T\|\bar x^rH\|_{L^2}.
\]
Density continuity and the weighted bound imply
\[
 \|\bar x^rH(t)\|_{L^2}
 \le \|H(t)\|_{L^2}^{1-r/a}
                 \|\bar x^aH(t)\|_{L^2}^{r/a}\longrightarrow0
                 \quad\text{as }t\downarrow0.
\]
The argument for \eqref{ws:eq:local-kinetic-tail} uses only the
strong-class bounds and prescribed traces, so both solutions satisfy
\[
 \sqrt\rho\,u(t),\
 \sqrt{\widetilde\rho}\,\widetilde u(t)
       \longrightarrow\sqrt{\rho_0}u_0
       \quad\hbox{strongly in }L^2\quad(t\downarrow0).
\]
Moreover,
\[
 \| (\sqrt\rho-\sqrt{\widetilde\rho})\widetilde u\|_{L^2}^2
 \le \int |H||\widetilde u|^2\,\dd
 \le C_T\|\bar x^rH(t)\|_{L^2}\longrightarrow0.
\]
Hence $\|\bar x^rH(0)\|_{L^2}=e(0)=0$, without prescribing $U$
on the initial vacuum set.

Subtracting the mass equations gives
\begin{equation}\label{ws:eq:uniqueness-density-equation}
 H_t+\widetilde u\cdot\nabla H
   +H\diver\widetilde u+\rho\diver U+U\cdot\nabla\rho=0.
\end{equation}
We claim that
\begin{equation}\label{ws:eq:uniqueness-density-bound}
 \|\bar x^rH(t)\|_{L^2}
 \le C_T\int_0^t\bigl(e(s)^{1/2}
                          +\|\nabla U(s)\|_{L^2}\bigr)\,\mathrm ds.
\end{equation}
Test \eqref{ws:eq:uniqueness-density-equation} by $H\bar x^{2r}$.
Since $|\nabla\log\bar x|\le C\bar x^{-1/2}$, transport contributes at most
\[
 C\bigl(\|\nabla\widetilde u\|_{L^\infty}
       +\|\widetilde u\bar x^{-1/2}\|_{L^\infty}\bigr)
                                      \|\bar x^rH\|_{L^2}^2.
\]
The term containing $\rho\diver U$ is bounded by
$C_T\|\bar x^rH\|_{L^2}\|\nabla U\|_{L^2}$.
Choose $p>2/(a-r)$ with $2p/(p-2)\le q$. Then
\[
\begin{aligned}
 \|\bar x^r U\cdot\nabla\rho\|_{L^2}
 &\le \|U\bar x^{-(a-r)}\|_{L^p}
       \|\bar x^a\nabla\rho\|_{L^{2p/(p-2)}}\\
 &\le C_T\bigl(e^{1/2}+\|\nabla U\|_{L^2}\bigr).
\end{aligned}
\]
Weighted Sobolev interpolation gives
\[
 \|\widetilde u\bar x^{-1/2}\|_{L^\infty}
 \le C_T\bigl(1+\|\nabla\widetilde u\|_{W^{1,q}}\bigr),
\]
with an $L^1(0,T)$ right-hand side. Gronwall, first applied to
$(\|\bar x^rH\|_{L^2}^2+\varepsilon)^{1/2}$, proves
\eqref{ws:eq:uniqueness-density-bound}. Setting
\[
 \mathcal I(t)=\int_0^t
             \bigl(e(s)+\|\nabla U(s)\|_{L^2}^2\bigr)\,\mathrm ds,
\]
yields
\begin{equation}\label{ws:eq:uniqueness-density-time-gain}
 \|\bar x^rH(t)\|_{L^2}^2\le C_T t\,\mathcal I(t).
\end{equation}

The momentum difference is
\begin{equation}\label{ws:eq:uniqueness-momentum-equation}
\begin{aligned}
 &\rho(U_t+u\cdot\nabla U)-\mu\Delta U
       -\nabla\bigl((\mu+\lambda(\rho))\diver U\bigr)\\
 &\quad=-\rho U\cdot\nabla\widetilde u
       -H(\widetilde u_t+\widetilde u\cdot\nabla\widetilde u)
       -\nabla\bigl(P(\rho)-P(\widetilde\rho)\bigr)\\
 &\qquad\quad
       +\nabla\bigl((\lambda(\rho)-\lambda(\widetilde\rho))
                                      \diver\widetilde u\bigr).
\end{aligned}
\end{equation}
On the bounded density range, local Lipschitz continuity gives
\begin{equation}\label{ws:eq:uniqueness-coefficient-differences}
\begin{aligned}
 &\|P(\rho)-P(\widetilde\rho)\|_{L^2}
       +\|\lambda(\rho)-\lambda(\widetilde\rho)\|_{L^2}\\
 &\qquad\le C_T\|H\|_{L^2}\le C_T\|\bar x^rH\|_{L^2}.
\end{aligned}
\end{equation}
Testing momentum by $U$, the coefficient differences contribute at most
\[
\begin{aligned}
 &C_T\|\bar x^rH\|_{L^2}
       \bigl(1+\|\nabla\widetilde u\|_{L^\infty}\bigr)
                                      \|\nabla U\|_{L^2}\\
 &\qquad\le\varepsilon\|\nabla U\|_{L^2}^2
       +C_{T,\varepsilon}\|\bar x^rH\|_{L^2}^2
                     \bigl(1+\|\nabla\widetilde u\|_{L^\infty}^2\bigr).
\end{aligned}
\]
Since $r/2>1/2$, \eqref{ws:eq:uniqueness-hardy} applies with exponent four:
\begin{align*}
 &\int |H||U|
      \bigl(|\widetilde u_t|
                  +|\widetilde u||\nabla\widetilde u|\bigr)\,\dd\\
 &\quad\le\|\bar x^rH\|_{L^2}\|U\bar x^{-r/2}\|_{L^4}\\
 &\qquad\quad{}\times
    \left(\|\widetilde u_t\bar x^{-r/2}\|_{L^4}
       +\|\nabla\widetilde u\|_{L^\infty}
                            \|\widetilde u\bar x^{-r/2}\|_{L^4}\right)\\
 &\quad\le\varepsilon\bigl(e+\|\nabla U\|_{L^2}^2\bigr)\\
 &\qquad\quad{}+C_{T,\varepsilon}\|\bar x^rH\|_{L^2}^2
        \left(\|\sqrt{\widetilde\rho}\,\widetilde u_t\|_{L^2}^2
          +\|\nabla\widetilde u_t\|_{L^2}^2
          +\|\nabla\widetilde u\|_{L^\infty}^2\right).
\end{align*}
The test is justified by $\chi_K DF_N(U)$ from
\eqref{ws:eq:physical-energy-balance}. For fixed $N$, the convective
boundary term is bounded by
\[
 CN\|\nabla\chi_K\|_{L^\infty}
                  \|\sqrt\rho u\|_{L^2}\,e^{1/2},
\]
and the stress and coefficient-difference boundary terms by
\[
 CN\|\nabla\chi_K\|_{L^2}
   \left(\|\nabla U\|_{L^2}
       +C_T\|\bar x^rH\|_{L^2}
                   \bigl(1+\|\nabla\widetilde u\|_{L^\infty}\bigr)\right).
\]
These vanish on $[\varepsilon,T]$ as $K\to\infty$.
Let $N\to\infty$, using $|DF_N(U)|\le|U|$ and the preceding interior bounds.
Using \eqref{ws:eq:uniqueness-density-time-gain} and absorbing gives
\begin{equation}\label{ws:eq:uniqueness-closed-difference}
\begin{aligned}
 e'+\mu\|\nabla U\|_{L^2}^2
 \le{}&C_T\bigl(1+\|\nabla\widetilde u\|_{L^\infty}\bigr)e\\
 &+C_T t\left(1+
        \|\sqrt{\widetilde\rho}\,\widetilde u_t\|_{L^2}^2
       +\|\nabla\widetilde u_t\|_{L^2}^2
       +\|\nabla\widetilde u\|_{L^\infty}^2\right)\mathcal I(t).
\end{aligned}
\end{equation}
By \eqref{ws:eq:strong-class},
\[
\begin{aligned}
 &\|\nabla\widetilde u\|_{L^\infty}\in L^1(0,T),\\
 &\int_0^T t\left(
   \|\sqrt{\widetilde\rho}\,\widetilde u_t\|_{L^2}^2
    +\|\nabla\widetilde u_t\|_{L^2}^2
    +\|\nabla\widetilde u\|_{L^\infty}^2\right)\,\mathrm dt<\infty,
\end{aligned}
\]
where the last term follows from
$\sqrt t\,\nabla\widetilde u\in L^2_tW^{1,q}_x\hookrightarrow L^2_tL^\infty_x$.
Thus the right-hand side of \eqref{ws:eq:uniqueness-closed-difference}
is integrable down to zero. Since $e(\varepsilon)\to0$, it follows that
\[
 \frac{d}{dt}\left(e+\frac\mu2\mathcal I\right)
 \le k_T(t)\left(e+\frac\mu2\mathcal I\right),
 \qquad k_T\in L^1(0,T).
\]
Gronwall gives $e=\nabla U=0$. Positive mass then gives $U=0$, and
\eqref{ws:eq:uniqueness-density-bound} gives $H=0$.

\begin{proof}[Proof of the well-posedness part of Theorem~\ref{ws:thm:main}]
Lemma~\ref{ws:lem:general-local} gives a local strong solution.
Proposition~\ref{ws:prop:apriori} gives \eqref{ws:eq:uniform-core}
and \eqref{ws:eq:flux-square-global} on its lifespan, and
Proposition~\ref{ws:prop:strong-continuation} extends it globally.
The preceding argument proves uniqueness. The large-time assertions
are proved in Subsection~\ref{ws:sec:asymptotics}.
\end{proof}

\section{Large-time behavior}\label{sec:large-time}

\subsection{Exponential relaxation on the torus}
Here $\Omega=\TT$.\label{per:sec:long-time}

We combine the pressure test from Subsection~\ref{per:sec:moment} with the
relative-energy argument of Fan, Li, and
Wang~\cite[Section~4.2]{FanLiWang2023}.

\begin{proposition}[Exponential relaxation]\label{per:prop:exponential-relaxation}
The global strong solution constructed in Subsection~\ref{per:sec:construction}
satisfies \eqref{per:eq:main-decay}. More precisely, for $2\le p<\infty$,
\begin{equation}\label{per:eq:smoothing-decay}
 \|\rho(t)-1\|_{L^{p}}+
 \min\{1,t\}^{1/2-1/p}\|u(t)-u_\infty\|_{W^{1,p}}
 \le C_pe^{-\eta_pt},\qquad t>0.
\end{equation}
For $1\le p<2$, the same bound holds without the time weight.
The constants depend only on $p$ and the low-order data and coefficient
bounds in Theorem~\ref{per:thm:main}.
\end{proposition}
\begin{proof}
Work in the zero-momentum frame of Subsection~\ref{per:sec:moment}.
The following tests are justified by time regularization.

For the relative energy
\[
 \mathcal E(t)=\int\left[
 \frac12\rho|u|^2+
 \frac{\rho^\gamma-1-\gamma(\rho-1)}{\gamma-1}\right]\dd,
\]
the energy identity, bounded density, and weighted Poincar\'e give
\begin{equation}\label{per:eq:relative-energy-equivalence}
 \mathcal E'+A_3^2=0,\qquad
 \mathcal E\asymp\|\sqrt\rho\,u\|_{L^{2}}^2+\|\rho-1\|_{L^{2}}^2
 \le C(A_3^2+\|\rho-1\|_{L^{2}}^2).
\end{equation}
To obtain density dissipation, apply \eqref{per:eq:pressure-test} and set
\[
 K(t)=\int\rho u\cdot\nabla\Delta^{-1}(\rho-1)\dd.
\]
As before, the projection and convection terms are bounded by $CA_3^2$.
Since the density is bounded,
\[
 \left|\int(2\mu+\lambda)\diver u\,(\rho-1)\dd\right|
 \le\varepsilon\|\rho-1\|_{L^{2}}^2+C_\varepsilon A_3^2.
\]
Since $(r^\gamma-1)(r-1)\ge(r-1)^2$ for $r\ge0$, we obtain
\[
 c\|\rho-1\|_{L^{2}}^2\le K'+CA_3^2,\qquad
 |K|\le C\|\sqrt\rho\,u\|_{L^{2}}\|\rho-1\|_{L^{2}}\le C\mathcal E.
\]
For a sufficiently small fixed $\delta>0$,
\[
 \mathcal E-\delta K\asymp\mathcal E,\qquad
 (\mathcal E-\delta K)'
       +c(A_3^2+\|\rho-1\|_{L^{2}}^2)\le0.
\]
Thus Gronwall and \eqref{per:eq:relative-energy-equivalence} give
\begin{equation}\label{per:eq:relative-energy-decay}
 \mathcal E(t)+\int_t^\infty A_3^2(s)\,\mathrm ds
 \le Ce^{-\eta t}
\end{equation}
for some $\eta>0$. Interpolation with the density upper bound yields
\begin{equation}\label{per:eq:density-decay}
 \|\rho(t)-1\|_{L^{p}}
 \le C_p e^{-\eta\min\{1/2,1/p\}t},
 \qquad 1\le p<\infty.
\end{equation}

It remains to obtain pointwise velocity decay. On the bounded
density range,
$\|P-\bar P\|_{L^{2}}^2+\|P-\bar P\|_{L^{4}}^4\le C\|\rho-1\|_{L^{2}}^2$.
The uniform bounds in \eqref{per:eq:main-bound} therefore reduce
\eqref{eq:log-energy-differential} to
\begin{equation}\label{per:eq:decaying-flux-differential}
 (A_1^2)'+cA_2^2\le C(A_3^2+\|\rho-1\|_{L^{2}}^2).
\end{equation}
By \eqref{per:eq:relative-energy-decay} and
\eqref{per:eq:energy-comparison},
$\int_t^\infty A_1^2\,\mathrm ds\le Ce^{-\eta t}$.
Average \eqref{per:eq:decaying-flux-differential} over starting times in
$[t-1,t]$. Integrating the resulting estimate from $t$ to infinity gives
\begin{equation}\label{per:eq:viscous-energy-decay}
 A_1^2(t)+A_3^2(t)+\int_t^\infty A_2^2(s)\,\mathrm ds
 \le Ce^{-\eta t},\qquad t\ge1.
\end{equation}
For $t\ge1$, the uniform bound for $A_2$ similarly reduces
\eqref{per:eq:acceleration-differential} to
\[
 (A_2^2)'+c\|\nabla\dot u\|_{L^{2}}^2
 \le C(A_3^2+\|\rho-1\|_{L^{2}}^2).
\]
The same averaging argument and \eqref{per:eq:viscous-energy-decay} yield
\begin{equation}\label{per:eq:acceleration-decay}
 A_2^2(t)+\int_t^\infty\|\nabla\dot u\|_{L^{2}}^2\,\mathrm ds
 \le Ce^{-\eta t},\qquad\text{for a.e. }t\ge2.
\end{equation}

For $2<p<\infty$, Hodge and Gagliardo--Nirenberg give the sharper bound
\[
 \|\nabla u\|_{L^{p}}
 \le C_p\bigl(A_1^{2/p}A_2^{1-2/p}
                  +A_3+\|\rho-1\|_{L^{p}}\bigr).
\]
Here the zero-mean parts of $G$ and $\omega$ are estimated using
\eqref{eq:flux-elliptic}, while $|\bar G|\le CA_3$.
Together with $\|u\|_{L^{p}}\le C_pA_3$, \eqref{per:eq:main-bound}, and
the exponential bounds above, this proves \eqref{per:eq:smoothing-decay}.
For $1\le p\le2$, use $\|u\|_{H^1}\le CA_3$ directly.
By positive-time $H^1$ continuity and weak lower semicontinuity in
$W^{1,p}$ for $p>2$, these estimates hold at every positive time.
Returning to the original frame gives \eqref{per:eq:main-decay}.
\end{proof}

\begin{remark}[Vacuum and density convergence]\label{per:rem:vacuum-long-time}
If $\rho_0(x_0)=0$, the transport formula gives
$\rho(t,X(t;x_0))=0$ at every finite time. Continuity of the density
then implies $\|\rho(t)-1\|_{L^{\infty}}\ge1$. Thus finite-$p$ exponential
convergence is compatible with persistent vacuum, whereas density
convergence in $L^\infty$ fails for these data.
\end{remark}

\subsection{Whole-space relaxation toward vacuum and algebraic decay}
\label{ws:sec:asymptotics}
Here $\Omega=\R^2$.

\begin{proof}[Proof of the large-time assertions in Theorem~\ref{ws:thm:main}]
We combine the uniform bounds \eqref{ws:eq:uniform-core} with
time-weighted pressure estimates, as in \cite[Lemma~3.4]{LiXin2019}.
In particular,
\[
 2\mu\le2\mu+\lambda(\rho)\le C,\qquad |\rho\lambda'(\rho)|\le C.
\]
The conservation laws established before
\eqref{ws:eq:physical-energy-balance} give
$M^{-1}\int\rho u\dd=u_\infty$.

\emph{Pressure moments and the first time weight.}
The pressure transport equation gives, for $p\ge2$,
\begin{equation}\label{ws:eq:decay-pressure-moment}
\begin{aligned}
 \frac{d}{dt}\|P\|_{L^p}^p
 &+(p\gamma-1)\int\frac{P^{p+1}}{2\mu+\lambda(\rho)}\dd\\
 &=-(p\gamma-1)\int\frac{P^pG}{2\mu+\lambda(\rho)}\dd .
\end{aligned}
\end{equation}
For $p=2,3$, Young's inequality and the elliptic interpolation bounds give
\[
\begin{aligned}
 \frac{d}{dt}\|P\|_{L^p}^p+c\|P\|_{L^{p+1}}^{p+1}
       &\le C A_1^2 A_2^{p-1},\\
 \|\nabla u\|_{L^{p+1}}^{p+1}
       &\le C A_1^2 A_2^{p-1}+C\|P\|_{L^{p+1}}^{p+1}.
\end{aligned}
\]
For a sufficiently large fixed $K_1$, set
\[
 F(t)=\int\left(
       \mu|\curl u|^2+(2\mu+\lambda(\rho))|\diver u|^2
                    -2P\diver u+K_1P^2\right)\dd .
\]
Completing the square and using Hodge give
\[
 c\bigl(\|\nabla u\|_{L^2}^2+\|P\|_{L^2}^2\bigr)
 \le F\le
 C\bigl(\|\nabla u\|_{L^2}^2+\|P\|_{L^2}^2\bigr),
 \qquad A_1^2\le CF.
\]
Testing momentum by $2\dot u$ and using continuity gives
\begin{equation}\label{ws:eq:corrected-energy-identity}
\begin{aligned}
 &\frac{d}{dt}\int\left(\mu\omega^2
       +(2\mu+\lambda(\rho))(\diver u)^2-2P\diver u\right)\dd
       +2A_2^2\\
 &\quad=-\mu\int\omega^2\diver u\dd
       +4\int G\det\nabla u\dd\\
 &\qquad-\int(2\mu+\lambda(\rho)+\rho\lambda'(\rho))(\diver u)^3\dd
       +2\gamma\int P(\diver u)^2\dd.
\end{aligned}
\end{equation}
Combining \eqref{ws:eq:corrected-energy-identity} with $K_1$ times
the pressure-square inequality and applying Young's inequality gives
\begin{equation}\label{ws:eq:decay-positive-energy}
 F'+c\bigl(A_2^2
                         +\|P\|_{L^3}^3\bigr)\le C F^2 .
\end{equation}
The global bounds already give $\int_1^\infty F\,dt\le C$.
Multiply \eqref{ws:eq:decay-positive-energy} by $t$.
The integrability of $F$ controls both the source and the coefficient in
$CF(tF)$. Gronwall therefore gives
\[
 \sup_{t\ge1}tF(t)+
 \int_1^\infty t\bigl(
       A_2^2+\|P\|_{L^3}^3\bigr)\,dt
 \le C.
\]

\emph{The second time weight.}
On the bounded density range, \eqref{eq:material-acceleration} gives
\[
 \frac{d}{dt}A_2^2
       +c\|\nabla\dot u\|_{L^2}^2
 \le C FA_2^2+C\|P\|_{L^4}^4.
\]
The right-hand side of \eqref{eq:material-acceleration} is bounded by
$C\int(|\nabla u|^4+P^2|\nabla u|^2)\dd$. Apply the fourth-power
estimate above.
Add a sufficiently large fixed $K_2$ times the pressure-cube
inequality. Then
\[
\begin{aligned}
 &\frac{d}{dt}\bigl(
       A_2^2+K_2\|P\|_{L^3}^3\bigr)
       +c\bigl(\|\nabla\dot u\|_{L^2}^2+\|P\|_{L^4}^4\bigr)\\
 &\qquad\le CF\bigl(
       A_2^2+K_2\|P\|_{L^3}^3\bigr).
\end{aligned}
\]
Multiply by $t^2$. The additional source is bounded by
$Ct(A_2^2+\|P\|_{L^3}^3)$, which is integrable by the preceding estimate.
Gronwall then yields
\begin{equation}\label{ws:eq:decay-time-weights}
\begin{aligned}
 &\sup_{t\ge1}t^2\bigl(
       A_2^2(t)
                              +\|P(t)\|_{L^3}^3\bigr)\\
 &\qquad+\int_1^\infty t^2\bigl(
       \|\nabla\dot u\|_{L^2}^2+\|P\|_{L^4}^4\bigr)\,dt\le C .
\end{aligned}
\end{equation}

\emph{Quantitative density decay.}
Lemma~\ref{ws:lem:potential} and $\log(1+s)\le s$ imply
\[
 \|G\|_{L^\infty}^2
 \le C\left(A_2^2
             +A_2
                                      \|\nabla\dot u\|_{L^2}\right).
\]
The two time-weighted estimates and Cauchy--Schwarz therefore give
\begin{equation}\label{ws:eq:decay-weighted-flux}
 \int_1^\infty t\|G(t)\|_{L^\infty}^2\,dt\le C,
 \qquad
 \int_t^\infty\|G(s)\|_{L^\infty}^2\,ds\le C/t
       \quad(t\ge1).
\end{equation}
Along every particle path, the pressure satisfies
\[
 \frac{dP}{dt}
   =-\frac{\gamma P(P+G)}{2\mu+\lambda(\rho)}
   \le-cP^2+C\|G(t)\|_{L^\infty}^2.
\]
Set
$f(t)=C\int_t^\infty\|G(s)\|_{L^\infty}^2\,ds$,
with the constant chosen as in this inequality.
Then $f(t)\le C/t$, and
\[
 \frac{d}{dt}(P+f)\le-cP^2
\]
along each path. Fix $t\ge2$.
If $P(t)+f(t)\le2f(t/2)$, then $P(t)\le C/t$.
Otherwise, throughout $[t/2,t]$, monotonicity of $P+f$ and of $f$
gives $P+f>2f(t/2)\ge2f$, hence $P\ge(P+f)/2$.
Consequently,
\[
 \frac{d}{ds}\frac1{P(s)+f(s)}\ge c/4,\qquad t/2\le s\le t,
\]
and integration gives $P(t)+f(t)\le8/(ct)$.
In both cases, $P(t)\le C/t$ with a constant independent of the path.
Vacuum paths satisfy $P=0$, and the flow covers all spatial points at
each finite time. Hence $\|\rho(t)\|_{L^\infty}\le Ct^{-1/\gamma}$.
For finite $p,r$, mass conservation gives
\[
 \|\rho(t)\|_{L^p}^p\le M\|\rho(t)\|_{L^\infty}^{p-1},\qquad
 \|P(t)\|_{L^r}^r
          \le M\|\rho(t)\|_{L^\infty}^{\gamma r-1}.
\]
Combined with the uniform bounds for $1\le t\le2$, these estimates give
\eqref{ws:eq:decay-density-pressure}.

\emph{Weighted Poincar\'e and gradient decay.}
The first time weight gives $\|\nabla u(t)\|_{L^2}^2\le C/t$.
Applying \eqref{ws:eq:mean-control-before-density} in the zero-momentum
frame and using the pressure $L^1$ rate,
we obtain
\[
 \frac12\int\rho|u-u_\infty|^2\dd+\frac1{\gamma-1}\int P\dd
 \le C\frac{\log(e+t)}t+C t^{-1+1/\gamma}
 \le C t^{-1+1/\gamma}.
\]
Here $\log(e+t)/t\le C_\gamma t^{-1+1/\gamma}$.
The energy balance gives $\int_t^S A_3^2\,ds\le Ct^{-1+1/\gamma}$ for
$S>t$. Letting $S\to\infty$ and using the pressure $L^2$ rate yields
\begin{equation}\label{ws:eq:decay-dissipation-tail}
 \int_t^\infty\bigl(
       \|\nabla u(s)\|_{L^2}^2+\|P(s)\|_{L^2}^2\bigr)\,ds
       \le C t^{-1+1/\gamma}.
\end{equation}
Since $F'\le CF^2$ and $\int_1^\infty F<\infty$,
Gronwall gives $F(t)\le CF(s)$ for $t/2\le s\le t$, $t\ge2$.
Averaging over this interval and using
\eqref{ws:eq:decay-dissipation-tail} gives
\[
 F(t)\le\frac Ct\int_{t/2}^t F(s)\,ds
             \le C t^{-2+1/\gamma}.
\]
Integrating \eqref{ws:eq:decay-positive-energy} from $t$ to infinity
then gives
\[
 \int_t^\infty\bigl(
       A_2^2+\|P\|_{L^3}^3\bigr)\,ds
       \le CF(t)+C\int_t^\infty F(s)^2\,ds
       \le C t^{-2+1/\gamma}.
\]
The differential inequality for
$A_2^2+K_2\|P\|_{L^3}^3$
has the same integrable coefficient $CF$.
The same averaging argument yields
\[
 A_2^2(t)
                  +K_2\|P(t)\|_{L^3}^3\le C t^{-3+1/\gamma}.
\]

Finally,
\[
\begin{aligned}
 \|G\|_{L^2}+\|\curl u\|_{L^2}
       &\le C t^{-1+1/(2\gamma)},\\
 \|\nabla G\|_{L^2}+\|\nabla\curl u\|_{L^2}
       &\le CA_2
        \le C t^{-3/2+1/(2\gamma)}.
\end{aligned}
\]
Interpolation gives, for every $2\le p<\infty$,
\[
 \|G\|_{L^p}+\|\curl u\|_{L^p}
       \le C_p t^{-3/2+1/(2\gamma)+1/p}.
\]
Since
$3/2-1/(2\gamma)-1/p\ge1-1/(\gamma p)$ for $p\ge2$,
the pressure rate and the div--curl estimate prove the gradient
bound in \eqref{ws:eq:decay-velocity}.

\emph{Removing the kinetic logarithm and improving acceleration.}
Introduce the natural expansion length and the centered moment radius
\[
 \ell(t)=(1+t)^{1/(2\gamma)},\qquad
 R_c(t)=\left[\frac1M\int\rho(x,t)
             (1+|x-u_\infty t|^2)^{a/2}\dd\right]^{1/a}\ge1.
\]

For any homogeneous $H^1$ field $v$ with finite weighted kinetic
energy, apply \eqref{ws:eq:log-weighted-poincare} after the change of
variables $x=u_\infty t+\ell(t)y$.
The rescaled density is $\ell(t)^2\rho(u_\infty t+\ell(t)y,t)$;
its mass is still $M$, and the Dirichlet integral of $v$ is unchanged
in two dimensions. Its $L^\gamma$ norm is at most
$M^{1/\gamma}(\ell^2\|\rho\|_{L^\infty})^{1-1/\gamma}
\le C(M+\ell^2\|\rho\|_{L^\infty})$. Consequently,
\[
\begin{aligned}
 \int\rho|v-v_\rho|^2\dd
 &\le C\left[\ell^2\|\rho\|_{L^\infty}+M
       +\int\rho\log\left(1+\frac{|x-u_\infty t|}{\ell}\right)\dd
             \right]\|\nabla v\|_{L^2}^2\\
 &\le C\left[\ell^2\|\rho\|_{L^\infty}
             +M\left(1+\log\left(1+\frac{R_c}{\ell}\right)\right)
             \right]\|\nabla v\|_{L^2}^2.
\end{aligned}
\]
The second inequality follows from Jensen's inequality for the
probability measure $\rho\,\mathrm dx/M$ and its moment bound of order $a$.
The constants are independent of the expansion length.

The continuity equation in the frame $x-u_\infty t$ gives
\[
\begin{aligned}
 aR_c^{a-1}|R_c'|
 &\le\frac aM\left(\int\rho|u-u_\infty|^2\dd\right)^{1/2}
       \left(\int\rho(1+|x-u_\infty t|^2)^{a-1}\dd\right)^{1/2}\\
 &\le\frac a{\sqrt M}R_c^{a-1}
                   \left(\int\rho|u-u_\infty|^2\dd\right)^{1/2}.
\end{aligned}
\]
The last step uses $0<2(a-1)/a<1$ and concavity with respect to the
same probability measure. Spatial truncation justifies differentiating
the moment. Thus $R_c$ is locally absolutely continuous and
\[
 |R_c'|\le M^{-1/2}
                \left(\int\rho|u-u_\infty|^2\dd\right)^{1/2}.
\]
Its value at time one is bounded by the initial moment and energy.

The established density rate gives
$\ell^2\|\rho\|_{L^\infty}\le C$, and conservation of momentum
gives $u_\rho=u_\infty$.
Applying the scaled inequality with $v=u$ and using
$\|\nabla u\|_{L^2}\le C\ell/(1+t)$ gives
\[
 (1+t)\frac{d}{dt}\left(\frac{R_c}{\ell}\right)
       +\frac1{2\gamma}\frac{R_c}{\ell}
 \le C\left[1+\log\left(1+\frac{R_c}{\ell}\right)\right]^{1/2}.
\]
The right-hand side grows more slowly than the linear term on the
left. Comparison with a sufficiently large constant thus proves
\begin{equation}\label{ws:eq:centered-moment-growth}
 \sup_{t\ge1}\frac{R_c(t)}{\ell(t)}\le C,\qquad
 \int\rho(1+|x-u_\infty t|^2)^{a/2}\dd
                       \le C(1+t)^{a/(2\gamma)}.
\end{equation}
In particular, the scaled inequality has a time-independent constant:
\begin{equation}\label{ws:eq:uniform-scaled-poincare}
 \int\rho|v-v_\rho|^2\dd\le C\|\nabla v\|_{L^2}^2,
 \qquad t\ge1.
\end{equation}
Taking $v=u$ removes the kinetic logarithm:
\[
 \int\rho|u-u_\infty|^2\dd
          \le C\|\nabla u\|_{L^2}^2
          \le C(1+t)^{-2+1/\gamma}.
\]

Apply \eqref{ws:eq:uniform-scaled-poincare} to the material acceleration.
The cutoff argument used for momentum conservation gives
\[
 \int\rho\dot u\dd=0.
\]
Indeed, $\rho\dot u\in L^1$, and the boundary term is bounded by
$C\|(G,\curl u)\mathbf1_{\{N\le|x|\le2N\}}\|_{L^2}\to0$.
Thus, for almost every $t\ge1$,
\begin{equation}\label{ws:eq:acceleration-poincare}
 A_2^2
                              \le C\|\nabla\dot u\|_{L^2}^2.
\end{equation}

Use this coercivity in the material-derivative inequality proved
above. The pressure $L^4$ rate gives
\[
\begin{aligned}
 \frac{d}{dt}A_2^2
             +c_0A_2^2
 &\le CF(t)A_2^2
                         +C\|P\|_{L^4}^4\\
 &\le CF(t)A_2^2
                         +C(1+t)^{-4+1/\gamma}.
\end{aligned}
\]
Since $F(t)\le C(1+t)^{-2+1/\gamma}$, this term can be absorbed
after a fixed time depending only on the data.
Integration from a fixed, sufficiently large $t_0$ gives
\[
\begin{aligned}
 A_2^2(t)
 &\le e^{-c(t-t_0)}
           A_2^2(t_0)\\
 &\quad+C\int_{t_0}^t
         e^{-c(t-s)}(1+s)^{-4+1/\gamma}\,ds\\
 &\le C(1+t)^{-4+1/\gamma}.
\end{aligned}
\]
For $t\ge2t_0$, split the integral at $t/2$: the first part is
exponentially small, while the second is bounded by
$C(1+t)^{-4+1/\gamma}$.
Consequently,
\begin{equation}\label{ws:eq:improved-acceleration-decay}
 A_2^2(t)
                         \le C(1+t)^{-4+1/\gamma},\qquad t\ge1.
\end{equation}
The uniform positive-time bounds cover the remaining finite interval.
These estimates hold for the standard positive-time energy representative,
with constants independent of the terminal time.

The estimates also give the qualitative density, gradient, and
acceleration limits. By momentum conservation,
\[
 \int\rho|u|^2\dd=M|u_\infty|^2+\int\rho|u-u_\infty|^2\dd
       \longrightarrow M|u_\infty|^2.
\]
\end{proof}

\begin{proposition}[Sharpness of the decay exponents]
\label{ws:prop:decay-sharpness}
The solution in Theorem~\ref{ws:thm:main} satisfies, for sufficiently large $t$,
\[
\begin{aligned}
 \|\rho(t)\|_{L^p}
 &\ge c_p(1+t)^{-\frac1\gamma(1-\frac1p)},
       &&1<p\le\infty,\\
 \|P(t)\|_{L^r}
 &\ge c_r(1+t)^{-1+\frac1{\gamma r}},
       &&1\le r\le\infty,\\
 \|\nabla u(t)\|_{L^p}
 &\ge c_p(1+t)^{-1+\frac1{\gamma p}},
       &&2\le p<\infty.
\end{aligned}
\]
These lower bounds have positive constants and match the upper bounds in
\eqref{ws:eq:decay-density-pressure}--\eqref{ws:eq:decay-velocity}. Moreover,
\[
\begin{aligned}
 0&<\limsup_{t\to\infty}(1+t)^{2-1/\gamma}
              \int\rho|u-u_\infty|^2\dd<\infty,\\
 0&<\limsup_{t\to\infty}(1+t)^{2-1/(2\gamma)}
              A_2(t)<\infty.
\end{aligned}
\]
The last two assertions give sharpness in the normalized-limsup
sense; no pointwise lower bound for these two quantities is asserted.
\end{proposition}
\begin{proof}
Use the expansion length $\ell$ and centered moment radius $R_c$ from the preceding proof.
The centered-moment estimate \eqref{ws:eq:centered-moment-growth} gives
$R_c(t)\le C(1+t)^{1/(2\gamma)}$.
At least half the mass lies in the ball centered at $u_\infty t$
of radius $2^{1/a}R_c(t)$. H\"older's inequality on this ball gives
the stated density lower bound for every finite $p>1$; for
$p=\infty$, use the mass divided by the area of the same ball.
Since $P=\rho^\gamma$, these bounds also give every pressure lower
bound. The corresponding upper bounds are
\eqref{ws:eq:decay-density-pressure}.

To obtain the gradient lower bound at exponent two, first apply
\eqref{ws:eq:uniform-scaled-poincare}.
For every compactly supported smooth vector test function $\varphi$,
subtract its density-weighted mean. Cauchy--Schwarz and that
Poincar\'e inequality give
\[
\begin{aligned}
 \left|\int\rho\dot u\cdot\varphi\dd\right|
 &=\left|\int\rho\dot u\cdot
       \left(\varphi-\frac1M\int\rho\varphi\dd\right)\dd\right|\\
 &\le CA_2\|\nabla\varphi\|_{L^2}.
\end{aligned}
\]
Thus
$\|\rho\dot u\|_{\dot H^{-1}}\le
CA_2$.
Equation~\eqref{ws:eq:improved-acceleration-decay}, the Fourier decomposition
$\rho\dot u=\nabla G+\mu\nabla^\perp(\curl u)$, and the elliptic
gradient estimate give
\[
\begin{aligned}
 \|G\|_{L^2}
 &\le C\|\rho\dot u\|_{\dot H^{-1}}
       \le C(1+t)^{-2+1/(2\gamma)},\\
 \|\nabla G\|_{L^2}
 &\le C\|\rho\|_{L^\infty}^{1/2}
                    A_2
       \le C(1+t)^{-2}.
\end{aligned}
\]
Interpolation yields
\[
 \|G(t)\|_{L^p}\le C_p(1+t)^{-2+1/(\gamma p)},
       \qquad 2\le p<\infty.
\]
Since $P=(2\mu+\lambda(\rho))\diver u-G$, the uniform coefficient
bound and the pressure lower bound give
\[
 c_p(1+t)^{-1+1/(\gamma p)}
 \le C\|\nabla u(t)\|_{L^p}
           +C_p(1+t)^{-2+1/(\gamma p)}.
\]
Absorb the last term for sufficiently large $t$ and use
\eqref{ws:eq:decay-velocity} to obtain the gradient assertion.

It remains to show that both limsups are positive. The density upper
bound and mass conservation imply
\[
 R_c(t)\ge c(1+t)^{1/(2\gamma)}.
\]
Indeed, a ball centered at $u_\infty t$ with radius a sufficiently small
fixed multiple of $(1+t)^{1/(2\gamma)}$ has mass at most $M/2$.
The mass outside this ball gives the stated moment lower bound. If the first normalized quantity tended to zero, the estimate for
$|R_c'|$ in the preceding proof would give
$|R_c'(t)|=o(t^{-1+1/(2\gamma)})$. Integration would give
$R_c(t)=o(t^{1/(2\gamma)})$, contradicting its lower bound.
The strengthened estimate \eqref{ws:eq:decay-kinetic} makes this
positive limsup finite.

Finally, continuity and the kinetic-energy cutoff argument used
for \eqref{ws:eq:physical-energy-balance} give
\[
 \frac{d}{dt}\int\rho|u-u_\infty|^2\dd
       =2\int\rho(u-u_\infty)\cdot\dot u\dd .
\]
Cauchy--Schwarz, first with a regularized square root, implies
\[
 \left|\frac{d}{dt}
       \left(\int\rho|u-u_\infty|^2\dd\right)^{1/2}\right|
       \le A_2.
\]
The relative kinetic energy tends to zero, and the improved
acceleration bound is integrable in time. Integration gives
\[
 \left(\int\rho(t)|u(t)-u_\infty|^2\dd\right)^{1/2}
 \le\int_t^\infty A_2(s)\,ds.
\]
If the second normalized quantity tended to zero, this integral
would be $o(t^{-1+1/(2\gamma)})$, contradicting the first limsup.
Its finiteness follows from
\eqref{ws:eq:improved-acceleration-decay}.
\end{proof}

\section*{Acknowledgments} 
The authors acknowledge the use of generative AI tools in developing and
checking mathematical arguments and in drafting and revising the manuscript.


\begin{thebibliography}{99}

\bibitem{CLMS1993}
R.~R.~Coifman, P.-L.~Lions, Y.~Meyer, and S.~Semmes,
Compensated compactness and Hardy spaces,
\emph{J. Math. Pures Appl.} (9) \textbf{72} (1993), 247--286.

\bibitem{CRW1976}
R.~R.~Coifman, R.~Rochberg, and G.~Weiss,
Factorization theorems for Hardy spaces in several variables,
\emph{Ann. of Math.} (2) \textbf{103} (1976), 611--635.
\href{https://doi.org/10.2307/1970954}
{doi:10.2307/1970954}.

\bibitem{FanLiWang2023}
X.~Fan, J.~Li, and X.~Wang,
Large-time behavior of the 2D compressible Navier--Stokes system in
bounded domains with large data and vacuum,
\href{https://arxiv.org/abs/2310.15520v1}{arXiv:2310.15520v1}, 2023.

\bibitem{HuangLi2016}
X.~Huang and J.~Li,
Existence and blowup behavior of global strong solutions to the
two-dimensional barotropic compressible Navier--Stokes system with
vacuum and large initial data,
\emph{J. Math. Pures Appl.} (9) \textbf{106} (2016), 123--154.
\href{https://doi.org/10.1016/j.matpur.2016.02.003}
{doi:10.1016/j.matpur.2016.02.003}.

\bibitem{HuangLi2022}
X.~Huang and J.~Li,
Global well-posedness of classical solutions to the Cauchy problem
of two-dimensional barotropic compressible Navier--Stokes system
with vacuum and large initial data,
\emph{SIAM J. Math. Anal.} \textbf{54} (2022), 3192--3214.
\href{https://doi.org/10.1137/21M1440943}
{doi:10.1137/21M1440943}.

\bibitem{JiuWangXin2014}
Q.~Jiu, Y.~Wang, and Z.~Xin,
Global well-posedness of 2D compressible Navier--Stokes equations
with large data and vacuum,
\emph{J. Math. Fluid Mech.} \textbf{16} (2014), 483--521.
\href{https://doi.org/10.1007/s00021-014-0171-8}
{doi:10.1007/s00021-014-0171-8}.

\bibitem{LiLiang2014}
J.~Li and Z.~Liang,
On local classical solutions to the Cauchy problem of the
two-dimensional barotropic compressible Navier--Stokes equations
with vacuum,
\emph{J. Math. Pures Appl.} \textbf{102} (2014), 640--671.
\href{https://doi.org/10.1016/j.matpur.2014.02.001}
{doi:10.1016/j.matpur.2014.02.001}.

\bibitem{LiXin2019}
J.~Li and Z.~Xin,
Global well-posedness and large time asymptotic behavior of
classical solutions to the compressible Navier--Stokes equations
with vacuum,
\emph{Ann. PDE} \textbf{5} (2019), Article~7.
\href{https://doi.org/10.1007/s40818-019-0064-5}
{doi:10.1007/s40818-019-0064-5}.

\bibitem{LiYang2026}
S.~Li and J.~Yang,
On global existence and large-time behaviour of weak solutions to the
compressible barotropic Navier--Stokes equations on $\mathbb T^2$ with
density-dependent bulk viscosity: beyond the Vaigant--Kazhikhov regime,
\href{https://arxiv.org/abs/2412.01226v4}{arXiv:2412.01226v4}, 2026.

\bibitem{Solonnikov1980}
V.~A.~Solonnikov,
Solvability of the initial-boundary-value problem for the equations
of motion of a viscous compressible fluid,
\emph{J. Soviet Math.} \textbf{14} (1980), 1120--1133.
\href{https://doi.org/10.1007/BF01562053}
{doi:10.1007/BF01562053}.

\bibitem{Trudinger1967}
N.~S.~Trudinger,
On imbeddings into Orlicz spaces and some applications,
\emph{J. Math. Mech.} \textbf{17} (1967), 473--483.
\href{https://doi.org/10.1512/iumj.1968.17.17028}
{doi:10.1512/iumj.1968.17.17028}.

\bibitem{VaigantKazhikhov1995}
V.~A.~Vaigant and A.~V.~Kazhikhov,
On existence of global solutions to the two-dimensional
Navier--Stokes equations for a compressible viscous fluid,
\emph{Siberian Math. J.} \textbf{36} (1995), 1108--1141.
\href{https://doi.org/10.1007/BF02106835}
{doi:10.1007/BF02106835}.

\end{thebibliography}
\end{document}